\documentclass[12pt]{imsart}

\usepackage[authoryear]{natbib}
\usepackage{ulem}

\RequirePackage[OT1]{fontenc}
\RequirePackage{amsthm,amsmath,amssymb}
\RequirePackage[colorlinks]{hyperref}
\usepackage[usenames]{color}
\usepackage{ulem}
\usepackage{memhfixc}  

\RequirePackage{amssymb}

\DeclareMathAlphabet{\mathonebb}{U}{bbold}{m}{n}
\newcommand{\1}{\ensuremath{\mathonebb{1}}}
\def\1{1\hskip-2.6pt{\rm l}}

\def\R{{\mathbb{R}}}

\def\E{{\mathbb{E}}}
\def\P{{\mathbb{P}}}
\def\IL{{\mathbb{L}}}

\def\SS{{\mathbb{S}}}
\def\FF{{\mathbb{F}}}

\def\rg{\mathrm{rg}}

\def\X{{\mathcal{X}}}

\def\I{{\mathcal I}}

\def\M{\mathcal M}

\def\KK{\mathcal K}
\def\HH{\mathcal H}

\def\EE{{\mathcal E}}
\def\LL{{\mathcal L}}
\def\RR{{\mathcal R}}

\def\NN{\mathcal N}

\def\expl{{\mathrm{ex}}}
\def\Cv{{\mathrm{Cv}}}
\def\crit{{\rm crit}}

\def\eps{{\varepsilon}}

\def\<{{\langle}}
\def\>{{\rangle}}

\newcommand{\pen}{\mathop{\rm pen}\nolimits}
\newcommand{\penn}{\mathop{\mathfrak{pen}}\nolimits}
\newcommand{\eref}[1]{(\ref{#1})}
\newcommand{\pa}[1]{\left({#1}\right)}
\newcommand{\norm}[1]{\left\|{#1}\right\|}
\newcommand{\cro}[1]{\left[{#1}\right]}
\newcommand{\ab}[1]{\left|{#1}\right|}
\newcommand{\ac}[1]{\left\{{#1}\right\}}

\newcommand{\bin}[2]{\ensuremath{\left(\!\!
\begin{array}{c}
#1\\
#2\\
\end{array}
\!\!\right)}}

\newtheorem{thm}{Theorem}
\newtheorem{lemma}{Lemma}
\newtheorem{prop}{Proposition}
\newtheorem{cor}{Corollary}

\newtheorem{ass}{Assumption}
\newtheorem{ex}{Example}

\def\EE{{\mathcal E}}

\def\I{{\mathcal I}}
\def\M{\mathcal M}

\def\ED{{\mathop{\rm pen}_{\Delta}\nolimits}}

\def\telque{\big |}

\begin{document}

\begin{frontmatter}
\title{Estimator selection in the Gaussian setting}
\runtitle{Estimator selection}

\begin{aug}
\author{\fnms{Yannick} \snm{Baraud}\ead[label=e1]{baraud@unice.fr}},
\author{\fnms{Christophe} \snm{Giraud}\ead[label=e2]{christophe.giraud@polytechnique.edu}}
\and
\author{\fnms{Sylvie} \snm{Huet}
\ead[label=e3]{sylvie.huet@jouy.inra.fr}
}

\runauthor{Y. Baraud et al}

\end{aug}



%
%
%
%
%

\begin{abstract}
We consider the problem of estimating the mean $f$ of a Gaussian vector $Y$ with independent components of common unknown variance $\sigma^{2}$. Our estimation procedure is based on estimator selection. More precisely, we start with an arbitrary and possibly infinite collection $\FF$ of estimators of $f$ based on $Y$ and, with the same data $Y$, aim at selecting an estimator among $\FF$ with the smallest Euclidean risk. No assumptions on the estimators are made and their dependencies with respect to $Y$ may be unknown. We establish a non-asymptotic risk bound for the selected estimator. 
As particular cases, our approach allows to handle the problems of aggregation and model selection as well as those of choosing  a window and a kernel for estimating a regression function, or tuning the parameter involved in a penalized criterion. We also derive oracle-type inequalities when $\FF$ consists of  linear estimators.  For illustration, we carry out two simulation studies. One aims at comparing our procedure to cross-validation for choosing a tuning parameter. The other shows how to implement our approach to solve the problem of variable selection in practice.
\end{abstract}

\end{frontmatter}
\maketitle
%
\section{Introduction}\label{sect:I}
\subsection{The setting and the approach}
We consider the Gaussian regression framework
$$Y_{i}=f_{i}+\eps_{i},\  i=1,\ldots,n$$
where $f=(f_{1},\ldots,f_{n})$ is an unknown vector of $\R^{n}$ and the $\eps_{i}$ are independent centered Gaussian random variables with common variance $\sigma^{2}$. Throughout the paper, $\sigma^2$ is assumed to be unknown which corresponds to the practical case. Our aim is to estimate $f$ from the observation of $Y$.
For specific forms of $f$, this setting allows to deal simultaneously with the following problems. 

\begin{ex}[Signal denoising]\label{ex-sd}
The vector $f$ is of the form
\[
f=(F(x_{1}),\ldots,F(x_{n}))  
\]
where $x_{1},\ldots,x_{n}$ are distinct points of a set $\X$ and $F$ is an unknown mapping from $\X$ into $\R$. 
\end{ex}

\begin{ex}[Linear regression]\label{ex-vs}
The vector $f$ is assumed to be of the form 
\begin{equation}\label{sv}
f=X\beta
\end{equation}
where 
$X$ is a $n\times p$ matrix, $\beta$ is an unknown $p$-dimensional
vector and $p$ some integer larger than 1 (and possibly larger than $n$). The columns of
the matrix $X$ are usually called predictors. When $p$ is large,
one may assume that the decomposition~\eref{sv} is sparse in the
sense that only few $\beta_{j}$ are non-zero. Estimating $f$ or finding
the predictors associated to the non-zero coordinates of
$\beta$ are classical issues. The latter is called variable selection.
\end{ex}

Our estimation strategy is based on estimator selection. More
precisely, we start with an arbitrary collection $\FF=\{\widehat f_{\lambda},\
\lambda\in\Lambda\}$ of estimators of $f$ based on $Y$ and aim at
selecting the one with the smallest Euclidean risk by using the same
observation $Y$. The way the estimators $\widehat f_{\lambda}$ depend on $Y$ may be arbitrary and possibly unknown. For example, the  $\widehat f_{\lambda}$ may be obtained from the minimization of a criterion, a Bayesian procedure or the guess of some experts.

\subsection{The motivation}
The problem of choosing some best estimator among a family of candidate ones is central in Statistics. Let us present some examples. 

\begin{ex}[Choosing a tuning parameter]\label{ex-ctp}
Many statistical procedures depend on a (possibly multi-dimensional) parameter $\lambda$ that needs to be tuned in view of obtaining an estimator with the best possible performance. For example, in the context of linear regression as described in Example~\ref{ex-vs},
the Lasso estimator (see Tibshirani~\citeyearpar{MR1379242} and Chen {\it et al.}~\citeyearpar{MR1639094}) defined by $\widehat f_{\lambda}=X\widehat{\beta_{\lambda}}$ with
\[
\widehat{\beta_{\lambda}}=\arg\!\min_{\beta\in\R^{p}}\cro{\norm{Y-X\beta}^{2}+\lambda \sum_{j=1}^{p}\ab{\beta_{j}}}
\]
depends on the choice of the parameter $\lambda\ge 0$. Selecting this parameter among a grid $\Lambda\subset \R_{+}$ amounts to selecting  a (suitable) estimator among the family $\FF=\{\widehat f_{\lambda},\ \lambda\in\Lambda\}$. 
\end{ex}

Another dilemma for Statisticians is the choice of a procedure to
solve a given problem. In the context of Example~\ref{ex-ctp}, there exist
many competitors to the Lasso estimator and one may alternatively
choose a procedure based on ridge regression (see Hoerl and
Kennard~\citeyearpar{ridge}), random forest or PLS 
(see Tenenhaus~\citeyearpar{MR1645125}, Helland~\citeyearpar{H01} and  Helland~\citeyearpar{H06}). 
Similarly, for the problem of signal denoising as described in Example~\ref{ex-sd}, popular approaches include spline smoothing, wavelet  decompositions and kernel estimators. The choice of a kernel may be possibly tricky. 

\begin{ex}[Choosing a kernel]
Consider the problem  described in Example~\ref{ex-sd} with $\X=\R$. For a kernel $K$ and a bandwidth $h>0$, the Nadaraya-Watson estimator (see Nadaraya~\citeyearpar{Nadaraya} and Watson~\citeyearpar{MR0185765}) $\widehat f_{K,h}\in\R^{n}$ is defined as
\[
\widehat f_{K,h}=\pa{\widehat F_{K,h}(x_{1}),\ldots,\widehat F_{K,h}(x_{n})}
\]
where for $x\in \R$
\[
\widehat F_{K,h}(x)={ \sum_{j=1}^n K\pa{x-x_{j}\over h} Y_{j}\over\sum_{j=1}^n K\pa{x-x_{j}\over h}}.
\]
There exist many possible choices for the kernel $K$, such as the Gaussian kernel
$K(x)=e^{-x^2/2}$, the uniform kernel $K(x)=\mathbf{1}_{|x|<1}$, etc. Given a (finite) family $\KK$ of candidate kernels $K$ and a grid $\HH\subset \R_{+}^{*}$ of possible values of $h$, one may consider the problem of selecting the best kernel estimator among the family $\FF=\{\widehat f_{\lambda},\ \lambda=(K,h)\in\KK\times\HH\}$. 
\end{ex}

\subsection{A look at the literature}
A common way to address the above issues is to use some cross-validation scheme such as leave-one-out or $V$-fold. Even though these resampling techniques are widely used in practice, little is known on their theoretical performances. For more details, we refer to Arlot and Celisse~\citeyearpar{arlot-2009} for a survey on cross-validation technics applied to model selection. Compared to these approaches, as we shall see, the procedure we propose is less time consuming and easier to implement. Moreover, it does not require to know how the estimators depend on the data $Y$ and we can therefore handle the following problem. 

\begin{ex}[Selecting among mute experts]
A Statistician is given a collection $\FF=\{\widehat f_{\lambda},\ \lambda\in\Lambda\}$ of estimators from a family $\Lambda$ of experts $\lambda$, each of which keeping secret the way his/her estimator $\widehat f_{\lambda}$ depends on the observation $Y$. The problem is then to find which expert $\lambda$ is the closest to the truth.  
\end{ex}

Given a selection rule among $\FF$, an important issue is to compare the risk of the selected estimator to those of the candidate ones. Results in this direction are available in the context of model selection, which can be seen as a particular case of estimator selection. More precisely, for the purpose of selecting a suitable model one starts with a collection $\SS$ of those, typically linear spaces chosen for their approximation properties with respect
to $f$, and one associates to each model $S\in\SS$ a suitable
estimator $\widehat f_{S}$ with values in $S$. Selecting a model then amounts to selecting an estimator among the collection $\FF=\{\widehat f_{S},\ S\in\SS\}$. 
For this problem, selection rules based on the minimization of a penalized criterion have been proposed in the regression setting by Yang~\citeyearpar{Yang99}, Baraud~\citeyearpar{MR1777129},  Birg\'e and Massart~\citeyearpar{MR1848946} and Baraud {\it et al}~\citeyearpar{BaGiHu2009}. Another way,  usually called Lepski's method, appears in a series of papers by Lepski~\citeyearpar{MR1091202,MR1147167,MR1214353,MR1191692}
and was originally designed to perform model selection among
collections of nested models. Finally, we mention that other
procedures based on resampling have interestingly emerged from the
work of Arlot~\citeyearpar{Arlot-These,MR2519533} and
C\'elisse~\citeyearpar{Celisse-these}. A common feature of those
approaches lies in the fact that the proposed selection rules apply to
specific collections of estimators only.

An alternative to {\it estimator selection} is {\it aggregation} which aims
at designing a suitable combination of given estimators in order to
outperform each of these separately (and even the best combination of
these) up to a remaining term.  Aggregation techniques  can be found
in Catoni~\citeyearpar{Catoni97,Catoni04}, Juditsky and Nemirovski~\citeyearpar{MR1792783},
Nemirovski~\citeyearpar{MR1775640}, Yang~\citeyearpar{MR1790617},
~\citeyearpar{MR1762904}, ~\citeyearpar{MR1946426},
Tsybakov~\citeyearpar{tsy03}, Wegkamp~\citeyearpar{Wegkamp03},
Birg\'e~\citeyearpar{MR2219712}, Rigollet and
Tsybakov~\citeyearpar{MR2356821}, Bunea, Tsybakov and
Wegkamp~\citeyearpar{MR2351101} and Goldenshluger~\citeyearpar{MR2488362}
for $\IL_{p}$-losses. Most of the aggregation procedures are based
on a sample splitting, one part of the data being used for building the estimators, the remaining part for selecting among these. Such a device requires that the observations be i.i.d. or at least that one has at disposal two independent copies of
the data. From this point of view our procedure differs from classical {\it
aggregation} procedures since we use the whole data $Y$ to build and select.  In the Gaussian regression setting that is  considered here, 
we mention the results of Leung
and Barron~\citeyearpar{MR2242356} for the problem of mixing
least-squares estimators. Their procedure uses the same data $Y$ to estimate and to aggregate but requires the variance to be known. Giraud~\citeyearpar{MR2543587} extends their results to the case where it is unknown.

\subsection{What is new here?}
Our approach for solving the problem of estimator selection is new. We introduce a collection $\SS$ of linear subspaces of $\R^{n}$ for approximating the estimators in $\FF$ and use a penalized criterion to compare them. As already mentioned and as we shall see, this approach requires no assumption on the family of estimators at hand and is easy to implement,
an {\tt R}-package being available on
\begin{center}
{\tt \small{\verb+http://w3.jouy.inra.fr/unites/miaj/public/perso/SylvieHuet_en.html+}}.
\end{center} 
A general way of comparing estimators in various statistical settings has been described in Baraud~\citeyearpar{hellinger}. However, the procedure proposed there is mainly abstract and inadequate in the Gaussian framework we consider.  

We prove a non-asymptotic risk bound for the estimator we select and
show that this bound is optimal in the sense that it essentially
cannot be improved (except for numerical constants maybe) by any other
selection rule. For the sakes of illustration and comparison, we apply
our procedure to various problems among which aggregation, model
selection, variable selection and selection among linear
estimators. In each of these cases, our approach allows to recover
classical results in the areas as well as to establish new ones. 
In
the context of aggregation we compute the
aggregation rates for the unknown variance case. These rates turn
out to be the same as those for the  known variance case.
For selecting an estimator among a family
of linear ones, we propose a new procedure and establish a risk bound  which
requires almost no assumption on the considered family.
 Finally, our approach
provides a way of selecting a suitable variable selection procedure
among a family of candidate ones.  It thus provides  an alternative to
cross-validation for which little is known.   

The paper is organized as follows. In Section~\ref{sect-main} we
present our selection rule and the theoretical properties of the
resulting estimator. For illustration, we show in Sections~\ref{sect-aggreg},~\ref{sect-linear} and~\ref{sect-vs} respectively, how the procedure can be used to aggregate preliminary estimators, select a linear estimator  
among a finite collection of candidate ones, or solve the problem of variable selection. Section~\ref{sec:numerique} is devoted to two simulation studies. One aims at comparing the performance of our procedure to the classical $V$-fold in view of selecting a tuning parameter among a grid. In the other, we evaluate the performance of the variable selection procedure we propose to some classical ones such as the Lasso, random forest, and others based on ridge and PLS regression.  Finally, the proofs are postponed to Section~\ref{sect-proof}.

Throughout the paper $C$ denotes a constant that may vary from line to line. 

\section{The procedure and the main result}\label{sect-main}
\subsection{The procedure}\label{desc-proc}
Given a collection 
$\FF=\{\widehat{f}_{\lambda}, \lambda \in
  \Lambda\}$ of estimators of $f$ based on $Y$,  the selection rule we propose is based on the choices of a family $\SS$ of linear subspaces of $\R^{n}$, a collection $\{\SS_{\lambda},\ \lambda\in\Lambda\}$ of (possibly random) subsets of $\SS$,   a weight function $\Delta$ and a penalty function $\pen$, both from $\SS$ into $\R_{+}$. We introduce those objects below and refer to Sections~\ref{sect-aggreg},~\ref{sect-linear} and~\ref{sect-vs}
  for examples. 

\subsubsection{The collection of estimators $\FF$}
The collection $\FF=\{\widehat{f}_{\lambda}, \lambda \in \Lambda\}$ can be arbitrary. In particular, $\FF$ need not be finite nor countable and it may consist of a mix of estimators based on the minimization of a criterion, a Bayes procedure or the guess of some experts. The dependency of these estimators with respect to $Y$ need not be known. Nevertheless, we shall see on examples how we can use this information, when available,  to improve the performance of our estimation procedure.

\subsubsection{The families $\SS$ and $\SS_{\lambda}$}
Let $\SS$ be a family of linear spaces of $\R^{n}$ satisfying the following. 
\begin{ass}\label{H0}
The family $\SS$ is finite or countable and for all $S\in\SS$, $\dim(S)\le n-2$.
\end{ass}
To each estimator $\widehat f_{\lambda}\in\FF$, we associate a (possibly random) subset $\SS_{\lambda}\subset\SS$. 

Typically, the family $\SS$ should be chosen to possess good approximation properties with respect to the elements of $\FF$ and $\SS_{\lambda}$ with respect to $\widehat f_{\lambda}$ specifically. One may take $\SS_{\lambda}=\SS$ but for computational reasons it will be convenient to allow $\SS_{\lambda}$ to be smaller. The choices of $\SS_{\lambda}$ may be made on the basis of the observation $\widehat f_{\lambda}$. We provide examples of $\SS$ and $\SS_{\lambda}$ in various statistical settings described in Sections~\ref{sect-aggreg} to~\ref{sect-vs}.

\subsubsection{The weight function $\Delta$ and the associated function $\ED$}
We consider a function $\Delta$ from $\SS$ into $\R_{+}$ and assume
\begin{ass}\label{H1}
\begin{equation}\label{sigma}
\Sigma=\sum_{S\in\SS}e^{-\Delta(S)}<+\infty.
\end{equation}
\end{ass}
Whenever $\SS$ is finite, inequality~\eref{sigma} automatically holds true. However, in practice $\Sigma$ should be kept to a reasonable size. When $\Sigma=1$, $e^{-\Delta(.)}$ can be interpreted as a prior distribution on $\SS$ and gives thus a Bayesian flavor to the procedure we propose. 
To the weight function $\Delta$, we associate the function $\ED$ mapping $\SS$ into $\R_{+}$ and defined by  
\begin{equation}
\E\cro{\pa{U-{\ED(S)\over n-\dim(S)}V}_{+}}=e^{-\Delta(S)}
\label{EDk.eq}
\end{equation}
where $x_{+}$ denotes the positive part of $x\in\R$ and $U,V$ are two independent $\chi^{2}$ random variables
with respectively $\dim(S)+1$ and $n-\dim(S)-1$ degrees of freedom. This function
can be easily computed from the quantiles of the Fisher distribution
as we shall see in Section~\ref{calcpen.st}. From a more theoretical point of view, it is shown in Baraud {\it et al}~\citeyearpar{BaGiHu2009} that under Assumption~\ref{H4} below, there exists a positive constant $C$ (depending on $\kappa$ only) such that
\begin{equation}\label{eq:oom}
\ED(S)\leq C(\dim(S)\vee \Delta(S)).
\end{equation}
\begin{ass}\label{H4}
There exists $\kappa\in(0,1)$ such that for all $S\in \SS$, 
\[
1\le \dim(S)\vee \Delta(S)\leq \kappa n.
\]
\end{ass}

\subsubsection{The selection criterion}
The selection procedure we propose involves a penalty function $\pen$ from $\SS$ into $\R_{+}$ with the following property.
\begin{ass}\label{H2}
The penalty function $\pen$ satisfies for some $K>1$,
\begin{equation}\label{pen}
\pen(S)\ge K\ED(S)\ \ \ {\rm for\ all}\ \ S\in\SS.
\end{equation}
\end{ass}
Whenever equality holds in~\eref{pen}, it derives from~\eref{eq:oom} that $\pen(S)$ measures the complexity of the model $S$ in terms of dimension and weight. 

Denoting $\Pi_{S}$ the projection operator onto a linear space $S\subset \R^{n}$, given the families $\SS_{\lambda}$, the penalty function $\pen$ and some positive number $\alpha$, we define 
\begin{equation}\label{eq:criterion}
\crit_{\alpha}(\widehat f_{\lambda})=\inf_{S\in\SS_{\lambda}}\cro{
\norm{Y-\Pi_{S}\widehat f_{\lambda}}^{2}+\alpha\norm{\widehat f_{\lambda}-\Pi_{S}\widehat f_{\lambda}}^{2}+\pen(S)\,\widehat\sigma^2_{S}},
\end{equation}
where
\begin{equation}\label{eq:var}
\widehat \sigma^2_{S}={\norm{Y-\Pi_{S}Y}^2\over n-\dim(S)}.
\end{equation}


\subsection{The main result}
For all $\lambda\in\Lambda$ let us set 
\begin{equation}\label{def-A}
A(\widehat f_{\lambda},\SS_{\lambda})=\inf_{S\in\SS_{\lambda}}\cro{\norm{\widehat f_{\lambda}-\Pi_{S}\widehat f_{\lambda}}^{2}+\pen(S)\,\widehat\sigma^2_{S}}.
\end{equation}
This quantity corresponds to an accuracy index for the estimator $\widehat f_{\lambda}$ with respect to the family $\SS_{\lambda}$. The following result holds. 
\begin{thm}\label{main}
Let $K>1,\alpha>0, \delta\ge 0$. Assume that Assumptions~\ref{H0},~\ref{H1} and~\ref{H2} hold. There exists a constant $C$ (given by~\eref{Cste.eq}) depending on $K$ and $\alpha$ only such that for any $\widehat f_{\widehat \lambda}$ in $\FF$  satisfying
\begin{equation}\label{def-est}
\crit_{\alpha}(\widehat f_{\widehat \lambda})\leq \inf_{\lambda\in\Lambda} \crit_{\alpha}(\widehat f_{\lambda})+\delta,
\end{equation}
we have the following bounds
\begin{eqnarray}
C\E\cro{\norm{f-\widehat f_{\widehat \lambda}}^{2}} &\le &
\E\cro{\inf_{\lambda\in\Lambda}\cro{\norm{f-\widehat
      f_{\lambda}}^{2}+A(\widehat
    f_{\lambda},\SS_{\lambda})}}+\Sigma\sigma^{2}+\delta \label{eqbase}\\
 &\le &\inf_{\lambda\in\Lambda}\ac{\E\cro{\norm{f-\widehat
       f_{\lambda}}^{2}}+\E\cro{A(\widehat
     f_{\lambda},\SS_{\lambda})}}+\Sigma\sigma^{2} +\delta \;\;\label{eq3}
\end{eqnarray}
 (provided that the quantity involved in the expectation
 in~\eref{eqbase} is measurable). 
Furthermore, if equality holds in~\eref{pen} and Assumption~\ref{H4} is satisfied, for each $\lambda\in\Lambda$
\begin{itemize}
\item if the set $\SS_{\lambda}$ is non-random, 
\end{itemize}
\vspace{-0.5cm}
\begin{eqnarray}
\lefteqn{C'\E\cro{A(\widehat f_{\lambda},\SS_{\lambda})}}\nonumber\\
&\le&\!\!\!\!\E\cro{\norm{f-\widehat f_{\lambda}}^{2}}+\inf_{S\in\SS_{\lambda}}\cro{\E\cro{\norm{\widehat f_{\lambda}-\Pi_{S}\widehat f_{\lambda}}^{2}}+(\dim(S)\vee \Delta(S))\sigma^{2}} \label{BorneA1}
\end{eqnarray}
\begin{itemize}
\item if there exists a (possibly random) linear space $\widehat S_{\lambda}\in\SS_{\lambda}$ such that $\widehat f_{\lambda}\in\widehat S_{\lambda}$ with probability 1, 
\end{itemize}
\begin{equation}\label{BorneA2}
C'\E\cro{A(\widehat f_{\lambda},\SS_{\lambda})}\le \E\cro{\norm{f-\widehat f_{\lambda}}^{2}}+\E\cro{\dim(\widehat S_{\lambda})\vee \Delta(\widehat S_{\lambda})}\sigma^{2},
\end{equation}
where $C'$ is a positive constant only depending on $\kappa$ and $K$.
\end{thm}

Let us now comment Theorem~\ref{main}. 

It turns out that inequality~\eref{eqbase} leaves no place for a substantial improvement in the sense that the bound we get is essentially optimal and cannot be improved (apart from constants) by any other selection rule among $\FF$. To see this, let us assume for simplicity that  $\FF$ is finite so that a measurable minimizer of $\crit_{\alpha}$ always exists and $\delta$ can be chosen as 0. Let $K=1.1$, $\alpha=1/2$ (to fix up the ideas), $\SS$ a family of linear spaces satisfying the assumptions of Theorem~\ref{main} and $\pen$, the penalty function achieving equality in~\eref{pen}. Besides, assume that $\SS$ contains a linear space $S$ such that $1\le \dim(S)\le n/2$ and associate to $S$ the weight $\Delta(S)=\dim(S)$. If $\SS_{\lambda}=\SS$ for all $\lambda$, we deduce from~\eref{eq:oom} and~\eref{eqbase} that for some universal constant $C'$, whatever $\FF$ and $f\in\R^{n}$ 
\begin{eqnarray}
\lefteqn{C'\E\cro{\norm{f-\widehat f_{\widehat \lambda}}^{2}}}\nonumber\\
&\le& \E\cro{\inf_{\lambda\in\Lambda}\cro{\norm{f-\widehat f_{\lambda}}^{2}+\inf_{S\in\SS}\pa{\norm{\widehat f_{\lambda}-\Pi_{S}\widehat f_{\lambda}}^{2}+\pen(S)\widehat \sigma^{2}_{S}}}}\nonumber\\
&\le& \E\cro{\inf_{\lambda\in\Lambda}\cro{\norm{f-\widehat f_{\lambda}}^{2}+\norm{\widehat f_{\lambda}-\Pi_{S}\widehat f_{\lambda}}^{2}+\dim(S)\widehat \sigma^{2}_{S}}}.\label{1000Borne}
\end{eqnarray}
In the opposite direction, the following result holds.
\begin{prop}\label{prop-opt}
There exists a universal constant $C$,  such that for any finite family $\FF=\{\widehat f_{\lambda},\ \lambda\in\Lambda\}$ of estimators and any selection rule $\widetilde \lambda$ based on $Y$ among $\Lambda$, there exists $f\in S$ such that
\begin{equation}\label{Borne1000} 
C\E\cro{\norm{f-\widehat f_{\widetilde \lambda}}^{2}}\ge \E\cro{\inf_{\lambda\in\Lambda}\cro{\norm{f-\widehat f_{\lambda}}^{2}+\norm{\widehat f_{\lambda}-\Pi_{S}\widehat f_{\lambda}}^{2}+\dim(S)\sigma^{2}}}.
\end{equation}
\end{prop}
We see that, up to the estimator $\widehat \sigma_{S}^{2}$ in place of $\sigma^{2}$ and numerical constants, the left-hand sides of~\eref{1000Borne} and~\eref{Borne1000} coincide.

In view of commenting~\eref{eq3} further, we continue assuming that $\FF$  is finite so that we can keep $\delta=0$ in~\eref{eq3}. A particular feature of~\eref{eq3} lies in the fact that the risk bound pays no price for considering a large collection $\FF$ of estimators. In fact, it is actually decreasing with respect to $\FF$ (or equivalently $\Lambda$) for the inclusion. This means that if one adds a new estimator to the collection $\FF$ (without changing neither $\SS$ nor the families $\SS_{\lambda}$ associated to the former estimators),  the risk bound for $\widehat f_{\widehat \lambda}$ can only be improved. In contrast, the computation of the estimator $\widehat f_{\widehat \lambda}$ is all the more difficult that $\ab{\FF}$ is large. More precisely, if the cardinalities of the families $\SS_{\lambda}$ are not too large, the computation of $\widehat f_{\widehat \lambda}$ requires around $\ab{\FF}$ steps.

The selection rule we use does not require to know how the estimators depend on $Y$. In fact, as we shall see, a more important piece of information is the ranges of the estimators $\widehat f_{\lambda}=\widehat f_{\lambda}(Y)$ as 
$Y$ varies in $\R^{n}$. A situation of special interest occurs when each $\widehat f_{\lambda}$ belongs to some (possibly random) linear space $\widehat S_{\lambda}$ in $\SS$ with probability one. By taking $\SS_{\lambda}$ such that $\widehat S_{\lambda}\in\SS_{\lambda}$ for all $\lambda$, we deduce from Theorem~\ref{main} by using~\eref{eq3} and~\eref{BorneA2} the following corollary.

\begin{cor}\label{maincor}
Assume that the Assumptions of Theorem~\ref{main} are satisfied, that Assumption~\ref{H4} holds and that equality holds in~\eref{pen}. If for all $\lambda\in\Lambda$ there exists a (possibly random) linear space $\widehat S_{\lambda}\in \SS_{\lambda}$ such that $\widehat f_{\lambda}\in\widehat S_{\lambda}$ with probability 1, then $\widehat f_{\widehat \lambda}$ satisfies 
\begin{equation}\label{eq21}
C\E\cro{\norm{f-\widehat f_{\widehat \lambda}}^{2}}\le \inf_{\lambda\in\Lambda}\cro{\E\cro{\norm{f-\widehat f_{\lambda}}^{2}}+\E\cro{\dim(\widehat S_{\lambda})\vee \Delta(\widehat S_{\lambda})}\sigma^{2}}+\delta,
\end{equation}
for some $C$ depending on $K$ and $\kappa$ only.
\end{cor}

One may apply this result in the context of model selection. One starts with a collection of models $\SS=\ac{S_{m},\ m\in\M}$ and associate to each  $S_{m}$ an estimator $\widehat f_{m}$ with values in $S_{m}$. By taking $\FF=\{\widehat f_{m},\ m\in\M\}$ (here $\Lambda=\M$) and $\SS_{m}=\ac{S_{m}}$ for all $m\in\M$, our selection procedure leads to an  estimator $\widehat f_{\widehat m}$ which satisfies
\begin{equation}\label{eq00}
C\E\cro{\norm{f-\widehat f_{\widehat m}}^{2}}\le \inf_{m\in\M}\cro{\E\cro{\norm{f-\widehat f_{m}}^{2}}+\pa{\dim(S_{m})\vee \Delta(S_{m})}\sigma^{2}}.
\end{equation}
When $\widehat f_{m}=\Pi_{S_{m}}Y$ for all $m\in\M$, our selection rule becomes 
\begin{equation}\label{critSM}
\widehat m=\arg\min_{m\in\M}\cro{\norm{Y-\widehat f_{m}}^{2}+\pen(S_{m})\,\widehat\sigma^2_{S_{m}}}
\end{equation}
and turns out to coincide with that described in Baraud {\it et al}~\citeyearpar{BaGiHu2009}. Interestingly, Corollary~\ref{maincor} shows that this selection rule can still be used for families $\FF$ of (non-linear) estimators of the form $\Pi_{S_{\widehat m}}Y$ where the $S_{\widehat m}$ are chosen randomly among $\SS$ on the basis of $Y$, doing thus as if the linear spaces $S_{\widehat m}$ were non-random. An estimator of the form $\Pi_{S_{\widehat m}}Y$ can be interpreted as resulting from a model selection procedures among the family of projection estimators $\{\Pi_{m}Y,\ m\in\M\}$ and hence,~\eref{critSM} can be used to choose some best model selection rule among a collection of candidate ones. 

\section{Aggregation}\label{sect-aggreg}
In this section, we consider the problems of {\it Model Selection Aggregation} (MS), {\it Convex Aggregation} ($\Cv$) and {\it Linear Aggregation} (L) defined below. Given $M\ge 2$ preliminary estimators of $f$, denoted $\ac{\phi_{k},\ k=1,\ldots,M}$, our aim is to build an estimator $\widehat f$ based on $Y$ whose risk is as close as possible to  $\inf_{g\in\FF_{\Lambda}}\norm{f-g}^{2}$ where 
\[
\FF_{\Lambda}=\ac{f_{\lambda}=\sum_{j=1}^{M}\lambda_{j}\phi_{j},\ \lambda\in\Lambda}
\]
and, according to the aggregation problem at hand,  $\Lambda$ is one of the three sets
\[
\Lambda_{{\rm MS}}=\ac{\lambda\in\{0,1\}^{M},\ \sum_{j=1}^{M}\lambda_{j}=1},\ \Lambda_{\Cv}=\ac{\lambda\in \R_{+}^{M},\ \sum_{j=1}^{M}\lambda_{j}=1},\ \Lambda_{{\rm L}}=\R^{M}.
\]
When $\Lambda=\Lambda_{{\rm MS}}$, $\FF_{\Lambda}$ is the set $\ac{\phi_{1},\ldots,\phi_{M}}$ consisting of the initial estimators. When $\Lambda=\Lambda_{\Cv}$, $\FF_{\Lambda}$ is the convex hull of the $\phi_{j}$. In the literature, one may also find 
\[
\Lambda_{\Cv}'=\ac{\lambda\in [0,1]^{M},\ \sum_{j=1}^{M}\lambda_{j}\le 1}
\]
in place of $\Lambda_{\Cv}$ in which case $\FF_{\Lambda}$ is the convex hull of $\ac{0,\phi_{1},\ldots,\phi_{M}}$. Finally, when $\Lambda=\Lambda_{{\rm L}}$, $\FF_{\Lambda}$ is the linear span of the $\phi_{j}$. 

Each of these three aggregation problems are solved {\it separately} if for each $\Lambda\in\ac{\Lambda_{{\rm MS}},\Lambda_{\Cv},\Lambda_{{\rm L}}}$ one can design an estimator $\widehat f=\widehat f(\Lambda)$ satisfying 
\begin{equation}\label{aggreg1}
\E\cro{\norm{f-\widehat f}^{2}}-C\inf_{g\in\FF_{\Lambda}}\norm{f-g}^{2}\le C'\psi_{n,\Lambda}\sigma^{2}
\end{equation}
with $C=1$, $C'>0$ free of $f,n,M$  and 
\begin{equation}\label{vaggreg}
\psi_{n,\Lambda}=\left\{
\begin{array}{cc}
M&  \mbox{if}\   \Lambda=\Lambda_{{\rm L}}\\
\sqrt{n\log(eM/\sqrt{n})} &   \mbox{if} \  \Lambda=\Lambda_{\Cv}\  \mbox{and}\ \sqrt{n}\le M\ \\
M & \mbox{if} \  \Lambda=\Lambda_{\Cv}\  \mbox{and}\ \sqrt{n}\ge  M \\
\log M &     \mbox{if} \  \Lambda=\Lambda_{{\rm MS}}.
\end{array}
\right.
\end{equation}
These problems have only been considered when the variance is known.
The quantity $\psi_{n,\Lambda}$ then corresponds to the best possible upper bound in~\eref{aggreg1} over all possible $f\in\R^{n}$ and preliminary estimators $\phi_{j}$ and is called the {\it optimal rate of aggregation}. For a more precise definition,  we refer the reader to Tsybakov~\citeyearpar{tsy03}. Bunea {\it et al}~\citeyearpar{MR2351101} considered the problem of solving these three problems {\it simultaneously} by building an estimator $\widehat f$ which satisfies~\eref{aggreg1} simultaneously for all $\Lambda\in\ac{\Lambda_{{\rm MS}},\Lambda_{\Cv},\Lambda_{{\rm L}}}$ and some constant $C>1$. This is an interesting issue since it is  impossible to know in practice which aggregation device should be used to achieve the smallest risk bound: as $\Lambda$ grows (for the inclusion), the bias $\inf_{g\in\FF_{\Lambda}}\norm{f-g}^{2}$ decreases while the rate $\psi_{n,\Lambda}$ increases.

The aim of this section is to show that our procedure provides a way of solving (or nearly solving) the three aggregation problems both {\it separately} and {\it simultaneously} when the variance is unknown.

Throughout this section, we consider the family $\overline \SS$ consisting of the $S_{m}$ defined for each $m\subset \ac{1,\ldots,M}$ and $m\neq \varnothing$ as the linear span of the $\phi_{j}$ for $j\in m$. Along this section, we shall use  the weight function $\Delta$ defined on $\overline \SS$ by
\[
\Delta(S_{m})=|m|+\log\cro{\bin{M}{|m|}}, 
\]
take $\alpha=1/2$ and $\pen(.)=1.1\ED(.)$ taking thus $K=1.1$. The choices of $\alpha$ and $K$ is only to fix up the ideas. Note that $\Delta$ satisfies Assumption~\ref{H1} with $\Sigma<1$. To avoid trivialities, we assume all along $n\ge 4$.  

\subsection{Solving the three aggregation problems separately}\label{AggregM}

\subsubsection{Linear Aggregation}
Problem (L) is the easiest to solve. Let us take $\FF=\FF_{\Lambda}$ with $\Lambda=\Lambda_{{\rm L}}$ and 
\begin{equation}\label{MLA}
\SS=\SS_{{\rm L}}=\ac{S_{\{1,\ldots,M\}}}
\end{equation} 
and $\SS_{\lambda}=\SS_{{\rm L}}$ for all $\lambda\in\Lambda_{{\rm L}}$. Minimizing $\crit_{\alpha}(f_{\lambda})$ over $f_{\lambda}\in\FF_{\Lambda}$ amounts to minimizing $\norm{Y-f_{\lambda}}^{2}$ over $f_{\lambda}\in S_{\{1,\ldots,M\}}$ and hence, the resulting estimator is merely $\widehat f_{{\rm L}}=\Pi_{S_{\{1,\ldots,M\}}}Y$. 
The risk of $\widehat f_{{\rm L}}$ satisfies
\[
\E\cro{\norm{f-\widehat f_{{\rm L}}}}\le \inf_{g\in \FF_{\Lambda}}\norm{f-g}^{2}+ M\sigma^{2}.
\]
whatever $n$ and $M$ which solves the problem of {\it Linear Aggregation}. 

\subsubsection{Model Selection Aggregation}
To tackle Problem (MS), we take $\FF=\FF_{\Lambda}$ with $\Lambda=\Lambda_{{\rm MS}}$, that is, $\FF_{\Lambda}=\ac{\phi_{1},\ldots,\phi_{M}}$, 
\begin{equation}\label{MMSA}
\SS=\SS_{{\rm MS}}=\{S_{\{1\}},\ldots,S_{\{M\}}\}
\end{equation}
and associate to each $f_{\lambda}=\phi_{j}$ the collection $\SS_{\lambda}$ reduced to $\ac{S_{\{j\}}}$. Note that $\dim(S)\le 1$ and $\Delta(S)=\log(eM)\ge \dim(S)$ for all $S\in\SS_{{\rm MS}}$, so that under the assumption that $\log(eM)\le n/2$ we may apply 
Corollary~\ref{maincor} with $\delta=0$ (since $\FF_{\Lambda}$ is finite), $\kappa=1/2$ and get that for some constant $C>0$ the resulting estimator $\widehat f_{{\rm MS}}$ satisfies
\[
C\E\cro{\norm{f-\widehat f_{{\rm MS}}}^{2}}\le \inf_{g\in \FF_{\Lambda}}\norm{f-g}^{2}+ \log(M)\sigma^{2}.
\]
This risk bound is of the form~\eref{aggreg1} except for the constant $C$ which is not equal to 1. We do not know whether Problem (MS) can be solved or not with $C=1$ when the variance $\sigma^{2}$ is unknown and $M$ is large (possibly larger than $n$). 

\subsubsection{Convex aggregation}\label{sect-CA}
For this problem, we emphasize the aggregation rate with respect to the quantity 
\begin{equation}\label{def-L}
L=\sup_{j=1,\ldots,M}{\norm{\phi_{j}}\over  \sigma\sqrt{n}}.
\end{equation}
If $M<\sqrt{n}L$, take again the estimator $\widehat f_{{\rm L}}$. Since the convex hull of the $\phi_{j}$ is a subset of the linear space $S_{\{1,\ldots,M\}}$,  for $\Lambda=\Lambda_{\Cv}$ we have
\[
\E\cro{\norm{f-\widehat f_{{\rm L}}}^{2}}\le \inf_{g\in \FF_{\Lambda}}\norm{f-g}^{2}+ M\sigma^{2}.
\]
Let us now turn to the case $M\ge\sqrt{n}L$. More precisely, assume that
\begin{equation}\label{HM}
2\le \sqrt{n}L\le M\le e^{-1}\min\ac{\sqrt{n}Le^{nL^{2}},e^{\sqrt{n}/(2L)}}
\end{equation}
and set $d(n,M)=n/(2\log(eM))$. We consider the family of estimators $\FF=\FF_{\Lambda}$ with $\Lambda=\Lambda_{\Cv}$ and
\begin{equation}\label{MCA}
\SS=\SS_{\Cv}=\SS_{\lambda}=\ac{S_{m}\in\overline \SS,\ |m|\le d(n,M)},\ \ \forall \lambda\in\Lambda_{\Cv}.
\end{equation}
The set $\Lambda_{\Cv}$ being compact, $\lambda\mapsto \crit_{\alpha}(f_{\lambda})$ admits a minimum $\widehat\lambda$ over $\Lambda_{\Cv}$ and we set $\widehat f_{\Cv}=\widehat f_{\widehat \lambda}$.
\begin{prop}\label{prop-convex}
There exists a universal constant $C>1$ such that
\[
\E\cro{\norm{f-\widehat f_{\Cv}}^{2}}-C\inf_{g\in\FF_{\Lambda}}\norm{f-g}^{2}\le C\sqrt{nL^{2}\log(eM/\sqrt{nL^{2}})}\sigma^{2}.
\]
\end{prop}
This risk bound is of the form~\eref{aggreg1} except for the constant $C$ which is not equal to 1. Again, we do not know whether Problem ($\Cv$) can be solved or not with $C=1$ when the variance $\sigma^{2}$ is unknown and $M$ possibly larger than $n$.

\subsection{Solving the three problems simultaneously}
Consider now three estimators $\widehat f_{{\rm L}},\widehat f_{{\rm MS}}, \widehat f_{\Cv}$ with values respectively in $S_{\{1,\ldots,M\}}$, $\bigcup_{j=1}^{M}S_{\{j\}}$ and the convex hull ${\mathcal C}$ of the $\phi_{j}$ (we use a new notation for this convex hull to avoid ambiguity). 
One may take the estimators defined in Section~\ref{AggregM} but any others would suit. The aim of this section is to select the one with the smallest risk to estimate $f$. To do so, we apply our selection procedure with $\FF=\{\widehat f_{{\rm L}},\widehat f_{{\rm MS}}, \widehat f_{\Cv}\}$, taking thus $\Lambda=\ac{{\rm L},{\rm MS},\Cv}$, and associate to each of these three estimators the families $\SS_{{\rm L}},\SS_{{\rm MS}},\SS_{\Cv}$ defined by~\eref{MLA},~\eref{MMSA} and~\eref{MCA} respectively and choose $\SS=\SS_{{\rm L}}\cup\SS_{{\rm MS}}\cup\SS_{\Cv}$. 
\begin{prop}\label{prop-simult}
Assume that~\eref{HM} holds and that $\log(eM)\le n/2$. There exists a universal constant $C>0$ such that whatever $\widehat f_{{\rm L}},\widehat f_{{\rm MS}}$ and $ \widehat f_{\Cv}$ with values in $S_{\{1,\ldots,M\}}$, $\bigcup_{j=1}^{M}S_{\{j\}}$ and ${\mathcal C}$ respectively, the selected estimator $\widehat f_{\widehat \lambda}$ satisfies for all $f\in\R^{n}$,
\[
C\E\cro{\norm{f-\widehat f_{\widehat \lambda}}^{2}}\le \inf_{\lambda\in \ac{{\rm L},{\rm MS},\Cv} }\cro{
\E\cro{\norm{f-\widehat f_{\lambda}}^{2}}
+B_{\lambda}},
\]
where 
\[
B_{{\rm L}}=\sigma^{2}M,\ \ B_{{\rm MS}}=\sigma^{2}\log M,\ \ B_{\Cv}= \sigma^{2}\cro{M\wedge\sqrt{nL^{2}\log(eM/\sqrt{nL^{2}})}}.
\]
In particular, if $\widehat f_{{\rm L}},\widehat f_{{\rm MS}}$ and $ \widehat f_{\Cv}$
fulfills~(\ref{aggreg1}), then  
\[
C\E\cro{\norm{f-\widehat f_{\widehat \lambda}}^{2}}\le \inf_{\lambda\in \ac{{\rm L},{\rm MS},\Cv} }\cro{
\inf_{g\in\FF_{\lambda}}\norm{f-g}^{2}+B_{\lambda}},
\]
where $\FF_{\lambda}$ stands for $\FF_{\Lambda}$ when $\Lambda=\Lambda_{\lambda}$.
\end{prop}
%

\section{Selecting among linear estimator}\label{sect-linear}
In this section, we consider the situation where  the estimators $\widehat f_{\lambda}$ are linear, that is, are  of the form $\widehat f_{\lambda}=A_\lambda Y$ for some known and deterministic $n\times n$ matrix $A_\lambda$. 
As mentioned before, this setting covers many popular estimation procedures including kernel ridge estimators, spline smoothing, Nadaraya estimators,  $\lambda$-nearest neighbors,  projection estimators,  low-pass filters, etc. In some cases $A_{\lambda}$ is symmetric (e.g. kernel ridge, spline smoothing, projection estimators), in some others $A_{\lambda}$ is non-symmetric and non-singular (as for Nadaraya estimators) and sometimes $A_{\lambda}$ can be both singular and non-symmetric (low pass filters, $\lambda$-nearest neighbors). A common feature of those procedures lies in the fact that they depend on a tuning parameter (possibly multidimensional) and their practical performances can be quite poor if this parameter is not suitably calibrated.
A series of papers have investigated the calibration of some of these procedures. To mention a few of them, Cao and Golubev~\citeyearpar{MR2301659} focus on spline smoothing, Zhang~\citeyearpar{MR2175849} on kernel ridge regression, Goldenshluger and Lepski~\citeyearpar{MR2449122} on kernel estimators  and   Arlot and Bach~\citeyearpar{arXiv:0909.1884}
propose a procedure to select among symmetric linear estimator with spectrum in $[0,1]$. The procedure we present can handle all these cases in an unified framework. 
Throughout the section, we assume that $\Lambda$ is finite.

\subsection{The families $\SS_{\lambda}$}
To apply our selection procedure, we need  to associate to each $A_{\lambda}$ a suitable collection of  approximation spaces $\SS_{\lambda}$. To do so, we introduce below a linear space $S_{\lambda}$ which plays a key role in our analysis.

For the sake of simplicity, let us first consider the case where $A_{\lambda}$ is non-singular. Then $S_{\lambda}$ is defined as the linear span 
of the right-singular vectors of $A^{-1}_{\lambda}-I$ associated to singular values smaller than 1. When $A_{\lambda}$ is symmetric, $S_{\lambda}$ is merely the linear span of the eigenvectors of $A_{\lambda}$ associated to eigenvalues not smaller than 1/2. If none of the singular values are smaller than 1, then $S_{\lambda}=\ac{0}$.

Let us now extend the definition of $S_{\lambda}$ to singular operators  $A_{\lambda}$. Let us recall that $\R^n=\ker(A_{\lambda})\oplus\rg(A_{\lambda}^*)$ where $A_{\lambda}^*$ stands for the transpose of $A_{\lambda}$ and $\rg(A_{\lambda}^*)$ for its range. The operator $A_{\lambda}$ then induces a one to one operator between $\rg(A^*_{\lambda})$ and $\rg(A_{\lambda})$. Write  $A_{\lambda}^+$ for the inverse of this operator from $\rg(A_{\lambda})$ to $\rg(A_{\lambda}^*)$. The orthogonal projection operator from $\R^{n}$ onto $\rg(A_{\lambda}^*)$ induces a linear operator from $\rg(A_{\lambda})$ into $\rg(A_{\lambda}^*)$, denoted $\overline \Pi_{\lambda}$. Then $S_{\lambda}$ is defined  as the linear span of the right-singular vectors of $A^+_{\lambda}-\overline \Pi_{\lambda}$ associated to singular values smaller than 1. Again if this set is empty,  $S_{\lambda}=\ac{0}$. When $A_{\lambda}$ is non-singular or symmetric, we recover the definition of $S_{\lambda}$ given above.

For each $\lambda\in\Lambda$, take $\SS_{\lambda}$ such that $\SS_{\lambda}\supset\ac{S_{\lambda}}$. From a theoretical point of view, it is enough to take $\SS_{\lambda}=\ac{S_{\lambda}}$ but practically it may be wise to use a larger set and by doing so, to possibly improve the approximation of $\widehat f_{\lambda}$ by elements of $\SS_{\lambda}$. One may for example take $\SS_{\lambda}=\ac{S_{\lambda}^1,\ldots,S_{\lambda}^{n-2}}$ where $S_{\lambda}^k$ is the linear span of the right-singular vectors associated to the $k$ smallest singular values of $A^+_{\lambda}-\overline \Pi_{\lambda}$.

\subsection{Choices of $\SS$, $\Delta$ and $\pen$}
Take $\SS=\bigcup_{\lambda\in\Lambda}\SS_{\lambda}$ and $\Delta$ of the form
\[
\Delta(S)=a\pa{1\vee\dim(S)}\ \ {\rm for\ all}\ S\in\SS
\]
where $a\ge 1$ satisfies Assumption~\ref{H1} with $\Sigma\le 1$. One may take $a=(\log \ab{\Lambda})\vee 1$ even though this choice is not necessarily the best. Finally, for some $K>1$, take $\pen(S)=K\ED(S)$ for all $S\in\SS$ and select $\widehat f_{\widehat \lambda}$ by minimizing the criterion given by~(\ref{eq:criterion}), taking thus $\delta=0$ in~\eref{def-est}.
%
%
%

\subsection{An oracle-type inequality for linear estimators}
The following holds. 
\begin{cor}\label{linear estimators}
Let $K>1$, $\kappa\in (0,1)$ and $\alpha>0$. If Assumption~\ref{H0} holds and $\Delta(S)\le \kappa n$ for all $S\in\SS$, the estimator $\widehat f_{\widehat \lambda}$ satisfies 
\[
Ca^{-1}\E\cro{\norm{f-\widehat f_{\widehat \lambda}}^{2}}\le \inf_{\lambda}\E\cro{\norm{f-\widehat f_{\lambda}}^{2}}+\sigma^{2},
\]
for some $C$ depending on $K,\alpha$ and $\kappa$ only. 
\end{cor}

The problem of selecting some best linear estimator among a family of those 
have also been considered in Arlot and Bach~\citeyearpar{arXiv:0909.1884} in the Gaussian regression framework, and in Goldenshluger and Lepski~\citeyearpar{MR2449122} in the multidimensional Gaussian white noise model. Arlot and Bach proposed a penalized procedure based on random penalties. Unlike ours, their approach requires that the operators be symmetric with eigenvalues in $[0,1]$ and  that the cardinality of $\Lambda$ is at most polynomial with respect to $n$. Goldenshluger and Lepski proposed a selection rule among families of kernel estimators to solve the problem of
structural adaptation. Their approach requires suitable assumptions on
the kernels while ours requires nothing. Nevertheless, we restrict to the case of the Euclidean loss whereas  Goldenshluger and Lepski considered more general $\IL_{p}$ ones.

\section{Variable selection}\label{sect-vs}
Throughout this section, we consider the problem of variable selection introduced in Example~\ref{ex-vs} and assume that $p\ge 2$ in order to avoid trivialities. When $p$ is small enough (say smaller than 20), this problem can be solved by using a suitable variable selection procedure that explores all the subsets of $\{1,\ldots,p\}$. For example, one may use the penalized criterion introduced in Birg\'e and Massart~\citeyearpar{MR1848946} when the variance is known, and the one in Baraud {\it et al}~\citeyearpar{BaGiHu2009} when it is not. When $p$ is larger, such an approach can no longer be applied since it becomes numerically intractable. To overcome this problem,  algorithms based on the minimization of convex criteria have been proposed among which are the Lasso, the Dantzig selector of Cand\`es and Tao~\citeyearpar{MR2382644}, the elastic net of Zou and Hastie~\citeyearpar{MR2137327}. An alternative to those criteria is the forward-backward algorithm described in Zhang~\citeyearpar{Zhang}, among others. Since there seems to be no evidence that one of these procedures outperforms all the others, it may be reasonable to mix them all and let the data decide which is the more appropriate to solve the problem at hand. As enlarging $\FF$  can only improve the risk bound of our estimator, only the CPU resources should limit the number of candidate estimators.

The procedure we propose could not only be used to select among those candidate procedures but also to select the tuning parameters they depend on. From this point of view, it provides an alternative to the cross-validation techniques which are quite popular but offer little theoretical guarantees. 

\subsection{Implementation roadmap}\label{sect-IRM}
Start by choosing a family $\LL$ of variable
selection procedures. Examples of such procedures are the Lasso, the Dantzig selector, the elastic net, among others. If necessary, associate to each $\ell\in\LL$ a
family  of
tuning parameters $H_{\ell}$. For example, in order to use the Lasso
procedure one needs to choose a tuning parameter $h>0$ among a grid
$H_{{\rm Lasso}}\subset \R_{+}$. If a selection procedure $\ell$
requires no choice of  tuning parameters, then one may take
$H_{\ell}=\ac{0}$. Let us denote by $\widehat m(\ell,h)$ the subset of
$\{1,\ldots,p\}$ corresponding to the predictors selected by the
procedure $\ell$ for the choice of the tuning parameter $h$.
For $m\subset \ac{1,\ldots,p}$, let $S_{m}$ be the linear span of the column vectors $X_{.,j}$ for $j\in m$ (with the convention $S_{\varnothing}=\ac{0}$). For $\ell \in\LL$ and $h\in H_{\ell}$, associate to the subset $\widehat m(\ell, h)$ an estimator $\widehat f_{(\ell,h)}$ of $f$ with values in $S_{\widehat m(\ell,h)}$ (one may for example take the projection of $Y$ onto the random linear space $S_{\widehat m(\ell,h)}$ but any other choice would suit). Finally, consider  the family $\FF=\{\widehat f_{\lambda},\
\lambda\in\Lambda\}$ of these estimators by taking
$\Lambda=\bigcup_{\ell\in\LL}(\{\ell\}\times H_{\ell})$ and set $\widehat \M=\ac{\widehat m(\lambda),\ \lambda\in\Lambda}$. All along we
assume that $\Lambda$ is finite 
(so that we take $\delta=0$ in~\eref{def-est}).

\subsubsection*{The approximation spaces and the weight function}
Throughout, we shall restrict ourselves to subsets of predictors with cardinality not larger than some $D_{\max}\le n-2$. In view of approximating the estimators $\widehat f_{\lambda}$, we suggest the collection $\SS$ given by 
\begin{equation}\label{Sm}
\SS=\bigcup \ac{S_{m}\telque\  m\subset \ac{1,\ldots,p}, {\rm card}(m)\leq D_{\max}}.
\end{equation}
We associate to $\SS$ the weight function $\Delta$ defined for $S\in\SS$ by
\begin{equation}
\Delta(S)=\log\cro{\bin{p}{D}}+\log(1+D)\ \ {\rm with}\ \ D=\dim(S).
\label{Delta.eq}
\end{equation}
Since 
\begin{eqnarray*}
\sum_{S\in\SS}e^{-\Delta(S)}&=&\sum_{D=0}^{p}\sum_{\scriptsize\begin{array}{cc} S\in\SS\\ \dim(S)=D\end{array}}e^{-\Delta(S)}\\
&\le& \sum_{D=0}^{p} e^{-\log(1+D)}\le 1+\log(1+p),
\end{eqnarray*}
Assumption~\ref{H1} is satisfied with $\Sigma=1+\log(1+p)$. 

Let us now turn to the choices of the $\SS_{\lambda}\subset \SS$. The criterion given by~\eref{eq:criterion} cannot be computed when $\SS_{\lambda}=\SS$ for all $\lambda$ as soon as $p$ is too large. In such a case, one must consider a smaller subset of $\SS$ and we suggest for $\lambda=(\ell,h)\in\Lambda$
\[
\SS_{(\ell,h)}=\ac{S_{\widehat m(\ell,h')},\ h'\in H_{\ell}}
\]
(where the $S_{m}$ are defined above), or preferably 
\[
\SS_{(\ell,h)}=\ac{S_{\widehat m(\ell',h')},\ \ell'\in \LL, h'\in H_{\ell}}
\]
whenever this latter family is not too large. Note that these two families are random. 


\subsection{The results}
Our choices of $\Delta$ and $\SS_{\lambda}$ ensure
that $\widehat f_{\lambda}\in S_{\widehat m(\lambda)}\in\SS_{\lambda}$
for all $\lambda\in\Lambda$ and that 
\[
\Delta(S_{\widehat m(\lambda)})\le 2\dim(S_{\widehat m(\lambda)})\log p.
\]
Hence, by applying Corollary~\ref{maincor} with $\widehat S_{\lambda}=S_{\widehat m(\lambda)}$, we get the following result.
\begin{cor}\label{cor-vs1}
Let $K>1$, $\kappa\in (0,1)$ and $D_{\max}$ be some positive integer satisfying $D_{\max}\leq \kappa n/ (2\log p)$. Let $\widehat \M=\{\widehat m(\lambda), \lambda\in\Lambda\}$ be a (finite) collection of random subsets of $\{1,\ldots,p\}$ with cardinality not larger than $D_{\max}$ based on the observation $Y$ and $\{\widehat f_{\lambda}, \lambda\in\Lambda\}$ a family of estimators $f$, also based on $Y$, such that $\widehat f_{\lambda}\in S_{\widehat m(\lambda)}$. By applying our selection procedure,  the resulting estimator $\widehat f_{\widehat \lambda}$ satisfies
\[
C\E\cro{\norm{f-\widehat f_{\widehat \lambda}}^{2}}\le \inf_{\lambda\in\Lambda}\cro{\E\cro{\norm{f-\widehat f_{\lambda}}^{2}}+\E\cro{\dim(S_{\widehat m(\lambda)})}\log(p)\sigma^{2}},
\]
where $C$ is a constant depending on the choices of $K$ and $\kappa$ only.
\end{cor}
Again, note that the risk bound we get is non-increasing with respect to $\Lambda$. This means that if one adds a new variable selection procedure or considers more tuning parameters to increase $\Lambda$, the risk bound we get can only be improved. 

Without additional information on the estimators $\widehat f_{\lambda}$ it is difficult to compare $\E\cro{\dim(S_{\widehat m(\lambda)})}\sigma^{2}$ and $\E\cro{\|f-\widehat f_{\lambda}\|^{2}}$. If $\widehat f_{\lambda}$ is of the form $\Pi_{S}Y$ for some deterministic subset $S\in\SS$ it is well-known that
\[
\E\cro{\norm{f-\Pi_{S}Y}^{2}}=\norm{f-\Pi_{S}f}^{2}+\dim(S)\sigma^{2}\ge \dim(S)\sigma^{2}.
\] 

Under the assumption that $f\in S_{m^{*}}$ and that $m^{*}$ belongs to $\widehat \M$ with probability close enough to 1, we can compare the risk of the estimator $\widehat f_{\widehat \lambda}$ to the cardinality of $m^*$.
  
\begin{cor}\label{cor-vs2}
Assume that the assumptions of Corollary~\ref{cor-vs1} hold and that $\widehat f_{\lambda}=\Pi_{S_{\widehat m(\lambda)}}Y$ for all $\lambda\in\Lambda$. If $f\in S_{m^{*}}$ for some non-void subset $m^{*}\subset \{1,\ldots,p\}$ with cardinality not larger than $D_{\max}$, then 
\[
C\E\cro{\norm{f-\widehat f_{\widehat \lambda}}^{2}}\le \log(p)|m^*|\sigma^2+R_{n}(m^{*})
\]
where $C$ is a constant depending on $K$ and $\kappa$ only, and 
\[
R_{n}(m^{*})= (\|f\|^2+n\sigma^2)\pa{\P\cro{m^*\not\in \widehat \M}}^{1/2}.
\]
\end{cor}
\cite{ZhaoYu06} gives sufficient conditions on the design $X$ to ensure that  $\P\cro{m^*\not\in \widehat \M}$ is exponentially small with respect to $n$ when the family $\widehat \M$ is obtained by using the LARS-Lasso algorithm with different values of the tuning parameter.

\section{Simulation study}\label{sec:numerique}
In the linear regression setting described in Example~\ref{ex-vs}, we carry out a simulation study  to evaluate the performances of our procedure to solve the two following problems. 

We first consider the problem, described in  Example~\ref{ex-ctp}, of  tuning the smoothing parameter of 
the Lasso procedure for estimating $f$. The performances of our procedure  are compared
with those of the $V$-fold cross-validation method. 
Secondly, we consider the problem of variable selection. We solve it by using our criterion in view of selecting among a  family $\LL$ of candidate variable selection procedures.  

Our simulation study is based on a large number of examples which have been chosen in view of covering a large variety of situations. Most of these have
been found in the literature in the context of Example~\ref{ex-vs} either for estimation or variable selection purposes when the number $p$ of predictors is
large. 

The section is organized as follows. The simulation design is  given in the following section. Then,  we describe how our procedure is applied for tuning the Lasso and performing variable selection. Finally,  we give the results of the
simulation study.

\subsection{\label{SimDes.st}Simulation design}

One example is determined by the number of observations $n$, the
number of variables $p$, the $n \times p$ matrix $X$, the values of
the parameters 
$\beta$, and the ratio signal/noise  $\rho$. It is denoted by 
$\expl(n,p,X,\beta, \rho)$, and the set of all considered examples is
denoted $\EE$. For each example, we carry out 400 simulations of $Y$
as a Gaussian random vector with expectation $f=X\beta$ and
variance $\sigma^{2} I_{n}$, where $I_{n}$ is the $n\times n$ identity
matrix, and $\sigma^{2} = \|f\|^{2}/n \rho$. 

The collection $\EE$ is composed of several collections $\EE_{e}$ for $e=1,\ldots,E$ where each collection $\EE_{e}$
is characterized by a vector of parameters
$\beta_{e}$, and a set $\X_{e}$ of matrices $X$:
\begin{equation*}
\EE_{e} = \left\{ \expl(n,p,X,\beta, \rho) : (n,p) \in \I, X \in
\X_{e}, \beta=\beta_{e}, \rho \in \RR \right\}
\end{equation*}
where $\RR= \{5, 10, 20\}$  and $\I$ consists of pairs $(n,p)$ such that
$p$ is smaller, equal or greater than $n$. 
The examples are described in further details  in
Section~\ref{SimEx.st}. They are inspired by examples found in
Tibshirani~\citeyearpar{MR1379242}, Zou and Hastie~\citeyearpar{MR2137327}, Zou~\citeyearpar{MR2279469}, and Huang
et al.~\citeyearpar{HMZ08} for comparing the Lasso method to the ridge, adaptive
Lasso and elastic net methods. They make up a large variety of
situations. 
They include cases where  
\begin{itemize}
\item the covariates are not, moderately or strongly  correlated, 
\item the covariates with zero
coefficients are weakly or highly
correlated with covariates with non-zero
coefficients, 
\item the covariates with non-zero
coefficients are grouped and correlated within these groups, 
\item  the
lasso method is known to be inconsistent,
\item  few  or many effects are present.
\end{itemize}

\subsection{\label{CVEDKhi.st} Tuning a smoothing parameter}
We consider here the problem of tuning the smoothing parameter of the Lasso estimator as described in Example~\ref{ex-ctp}. Instead of considering the Lasso estimators for a fixed grid $\Lambda$ of smoothing parameters $\lambda$, we rather focus on the sequence $\{\widehat f_1,\ldots,\widehat f_{D_{\max}}\}$ of estimators given by the $D_{\max}$ first steps of the LARS-Lasso algorithm proposed by Efron {\it et al.}~\citeyearpar{MR2060166}. Hence, the tuning parameter is here the number $h\in H=\ac{1,\ldots,D_{\max}}$ of steps. In our simulation study, we compare the performance of our criterion to that of the $V$-fold cross-validation for the problem of selecting the best estimator among the collection $\FF=\{\widehat f_1,\ldots,\widehat f_{D_{\max}}\}$.

\subsubsection{The estimator of $f$ based on our procedure}
We recall that our selection procedure relies on the choices of families $\SS$, $\SS_{h}$ for $h\in H$, a weight function $\Delta$, a penalty function $\pen$ and two universal constants $K>1$ and $\alpha>0$. We choose the family $\SS$ defined by~\eref{Sm}.
We associate to $\widehat f_{h}$ the family $\SS_{h}=\{S_{\widehat m(h')}|\ h'\in H\}\subset \SS$ where the $S_m$ are defined in Section~\ref{sect-IRM} and $\widehat m(h')\subset \ac{1,\ldots,p}$ is the set of indices corresponding to the predictors retuned by the LARS-Lasso algorithm at step $h'\in H$. We take $\pen(S) = K \ED(S) $ with $\Delta(S)$ defined by~\eref{Delta.eq} and $K=1.1$. This value of $K$ is consistent with what is suggested in Baraud {\it et al.}~\citeyearpar{BaGiHu2009}. The choice of $\alpha$ is based on the following considerations. First, choosing $\alpha$ around one seems
reasonable since  it weights 
similarly  the term $\|Y-\Pi_{S}\widehat{f}_{\lambda}\|^{2}$ which measures how well 
the estimator fits the data and the approximation term $\|\widehat{f}_{\lambda}-\Pi_{S}\widehat{f}_{\lambda}\|^{2}$ involved in our criterion~\eref{eq:criterion}. Second,  simple calculation shows that
the constant $C^{-1}=C^{-1}(1.1, \alpha)$ involved in Theorem~\ref{main} is minimum for
$\alpha$ close to 0.6. We therefore carried out our
simulations for $\alpha$ varying from 0.2 to 1.5. The
results being very similar for $\alpha$ between 0.5 and 1.2, we  choose $\alpha=0.5$. We denote by $\widehat{f}_{\ED}$ the resulting estimator of $f$.

\subsubsection{The estimator of $f$ based on $V$-fold cross-validation}
For each $h \in H$, the prediction error is estimated using a $V$-fold
cross-validation procedure, with $V=n/10$. The estimator
$\widehat{f}_{CV}$  is chosen by minimizing the estimated prediction error.

\subsubsection{The results}
The simulations were carried out with {\tt R} ({\tt www.r-project.org}) using the library
{\tt elasticnet}.

For each example $\expl\in\EE$, we estimate on the basis of 400 simulations the oracle risk  
\begin{equation}
O_{\expl} = \E\left( \min_{h \in H} \|f - \widehat{f}_{h}\|^{2}\right),
\label{Oracl.eq}
\end{equation}
and the Euclidean risks $R_{\expl}(\widehat{f}_{\ED})$ and $R_{\expl}(\widehat{f}_{CV})$ of  $\widehat{f}_{\ED}$ and $\widehat{f}_{CV}$ respectively.

The 
results presented in Table~\ref{tb1} show that our procedure tends to
choose a better estimator than the CV  in the sense that the
ratios $R_{\expl}(\widehat{f}_{\ED})/O_{\expl}$ are closer to one than $R_{\expl}(\widehat{f}_{CV})/O_{\expl}$.

\begin{table}
\begin{tabular}{r|rrrrrrr}
&&&\multicolumn{5}{c}{quantiles} \\
procedure& mean & std-err & $0\%$ & $50\%$ & $75\%$ & 99\% & $100\%$ \\
CV & 1.18 & 0.08 & 1.05 & 1.18 & 1.24 & 1.36 & 1.38 \\
$\ED$  & 1.065 & 0.06 & 1.01  &  1.055 & 1.084 & 1.18 & 2.27\\
\end{tabular}
\caption{\label{tb1} Mean, standard-error and quantiles of the ratios
$R_{\expl}/O_{\expl}$ calculated over 
all $\expl \in \EE$ such that $O_{\expl} < n \sigma^{2} /3$. The
number of such examples equals 654, see Section~\ref{SimEx.st}.
}
\end{table}


Nevertheless, for a few
examples these ratios are larger for our
procedure than for the CV. 
These examples correspond to situations where the Lasso estimators are highly biased.

In practice, it is worth considering several estimation procedures in order to increase the chance to have good estimators of $f$ among the family $\FF$. Selecting among candidate procedures is the purpose of the following simulation experiment in the variable selection context.

\subsection{\label{SelVar.st}Variable selection}
In this section, we consider the problem of variable selection and use the procedure and notations introduced in Section~\ref{sect-vs}. To solve this problem, we consider estimators of the form $\widehat f_{\widehat m}=\Pi_{S_{\widehat m}}Y$ where $\widehat m$ is a random subset of $\{1,\ldots,p\}$ depending on $Y$. Given a family 
$\widehat\M=\{\widehat m(\ell,h),\ \widehat m(\ell,h) \in \LL\times H_{\ell}\}$ of such random sets, we consider the family $\FF=\{\widehat f_{\widehat m(\ell,h)}|\ (\ell,h)\in\LL\times H_{\ell}\}$. The descriptions of $\LL$ and $H_{\ell}$ are postponed to Section~\ref{proc.st}. Let us merely mention that we choose $\LL$ which gathers variable selection procedures based on the Lasso, ridge regression, Elastic net, PLS1 regression, Adaptive Lasso, Random Forest, and on an exhaustive research among the subsets of $\ac{1,\ldots,p}$ with small cardinality. For each procedure $\ell$, the parameter set $H_{\ell}$ corresponds to different choices of tuning parameters. 
%
%
%
For each $\lambda=(\ell,h)\in \LL\times H_{\ell}$, we take $\SS_{\lambda}=\{S_{\widehat m(\ell,h)}\}$ so that our selection rule over $\FF$ amounts to minimizing over $\widehat\M$
\begin{equation}\label{critP.eq}
\crit(m)=\|Y-\Pi_{S_{m}}Y\|^2+K\pen_{\Delta}(S_m)\widehat\sigma_{S_m}^2,
\end{equation}
where $\pen_{\Delta}$ is given by~\eref{EDk.eq}. 

\subsubsection{Results}
The simulations were carried out  with {\tt R} ({\tt
www.r-project.org}) using the libraries {\tt elasticnet}, {\tt
randomForest}, {\tt pls} and  the program {\tt lm.ridge} in the library
{\tt MASS}. We first select the tuning parameters associated to the procedures $\ell$ in $\LL$. More precisely, for each $\ell$ we  select an estimator among the collection $\FF_{\ell}=\{\widehat f_{\widehat m(\ell,h)}|\ h\in H_{\ell}\}$ by minimizing Criterion~\eref{critP.eq} over $\widehat\M_{\ell}=\left\{{\widehat m(\ell,h)}| h \in H_{\ell}\right\}$.
We denote by $\widehat m(\ell)$ the selected set and by $\widehat f_{\widehat m(\ell)}$ the corresponding projection estimator. For each example $\expl \in \EE$ and each method $\ell \in \LL$, we  estimate the risk 
\begin{equation*}
R_{\expl, \ell} = \E\left(  \|f - \widehat f_{\widehat m(\ell)}\|^{2}\right)
\end{equation*}
of $\widehat f_{\widehat m(\ell)}$ on the basis of 400 simulations and we do the same  
to calculate that of our estimator $\widehat f_{\widehat m}$,
\begin{equation*}
R_{\expl, \mbox{all}} = \E\left(  \|f - \widehat f_{\widehat m} \|^{2}\right).
\end{equation*}
Let us now define the minimum of these risks over all methods: 
\begin{equation*}
R_{\expl, \min} = \min\left\{ R_{\expl,\mathrm{all}},  R_{\expl,
   \ell}, \ell \in \LL\right\}.
\end{equation*}

We compare the ratios $R_{\expl,\ell}/R_{\expl, \min}$ 
for $\ell\in\LL\cup\{{\rm all}\}$ to judge the performances of the candidate procedures on each example $\expl\in\EE$. The mean, standard deviations and quantiles of the sequence $\{R_{\expl,\ell}/R_{\expl, \min},\ \expl\in\EE\}$ are presented  in Table~\ref{tb2}. In particular, the results show that 
\begin{itemize}
\item none of the procedures $\ell$ in $\LL$ outperforms all the others simultaneously over all examples,
\item our procedure, corresponding to $\ell={\rm all}$, achieves the  smallest mean value. Besides, this value is very close to one.
\item the variability of our procedure is small compared to the others
\item for all examples, our procedure selects an estimator the risk of which does not exceed twice that of the oracle.
\end{itemize}

\begin{table}
\begin{tabular}{l|rrrrrr}
&&&\multicolumn{4}{c}{quantiles} \\ \cline{4-7}
method & mean & std-err &  $50\%$ & $75\%$ & 95\% & $100\%$ \\ \hline
Lasso  &  2.82&  9.40 &     1.12 &   1.33  & 6.38 &127\\
ridge  &  1.76&  1.90 &     1.42 &   1.82  & 2.87 & 36.9\\
pls    &  1.50&  1.20 &     1.22 &   1.50  & 2.58 & 17\\
en     &  1.46&  1.90 &     1.12 &   1.33  & 2.57 & 29\\
ALridge&  1.20&  0.31 &     1.15 &   1.26  & 1.51 &  5.78\\
ALpls  &  1.29&  0.87 &     1.14 &   1.29  & 1.75 & 12.7\\
rFmse  &  4.13&  9.50 &     1.38 &   2.04  & 19.2 &118\\
rFpurity& 3.99& 10.00 &     1.42 &   2.06  & 15.1 & 138\\
exhaustive & 22.9 & 45 &    6.30 &    24.5 & 92.9 & 430 \\
all  &    1.16&  0.16 &     1.12 &   1.25  &  1.47 & 1.95\\

\end{tabular}
\caption{\label{tb2} For each $\ell\in\LL\cup\ac{{\rm all}}$, mean,
standard-error and quantiles of the ratios $R_{\expl, \ell} /R_{\expl, \min}$
 calculated over
all $\expl \in \EE$. The number of examples in the collection $\EE$ is equal to 660. 
}
\end{table}

The false discovery rate (FDR)
and the true discovery rate (TDR) are also parameters of interest in the context of variable selection. These quantities 
are given at Table~\ref{tb3} for each example when $\rho=10$ and
$n=p=100$. Except for one example, the FDR is
small, while the TDR is varying a lot among the examples. 

\begin{table}
\begin{tabular}{l|rrrrrrrrrrr}
~&$\EE_{1}$&$\EE_{2}$&$\EE_{3}$&$\EE_{4}$&$\EE_{5}$&$\EE_{6}$&$\EE_{7}$&$\EE_{8}$&$\EE_{9}$&$\EE_{10}$&$\EE_{11}$ \\
FDR ~&0.045    &0.026    &0.004    &0.026    &0.018    &0.041    &0.012
&0.026&0.042& 0.15&0.014\\
TDR &0.74     &0.63     &0.18     &0.63     &0.17     &0.99     &1
&1&0.98&0.29&0.20\\
\end{tabular}
\caption{\label{tb3} ~ False dicovery rate (FDR) and true discovery
rate (TDR) using our method, for each example with $\rho=10$ and
$n=p=100$.}
\end{table}





\

\section{Proofs}\label{sect-proof}
\subsection{Proof of Theorem~\ref{main}}
Throughout this section, we use the following notations. For all $\lambda\in\Lambda$ and  $S\in\SS_{\lambda}$, we write
\[
\crit_{\alpha}(\widehat f_{\lambda},S)=\norm{Y-\Pi_{S}\widehat f_{\lambda}}^{2}+\sigma^2\penn(S)+\alpha\norm{\widehat f_{\lambda}-\Pi_{S}\widehat f_{\lambda}}^{2},
\]
where 
\begin{equation}\label{penn}
\penn(S)=\pen(S)\,\widehat\sigma^2_{S}/\sigma^2,\ \ \textrm{for all }S\in\SS.
\end{equation}
For all $\lambda\in\Lambda$, let $S(\lambda)\in \SS_{\lambda}$ be such that
\[
\crit_{\alpha}(\widehat f_{\lambda},S(\lambda))\leq \crit_{\alpha}(\widehat f_{\lambda})+\delta.
\]
We also write $\eps=Y-f$ and $\overline S$ for the linear space generated by $S$ and $f$.
It follows the facts that for all $\lambda\in\Lambda$ and $S\in\SS_{\lambda}$
\[
\crit_{\alpha}(\widehat f_{\widehat \lambda},S(\widehat\lambda))\le \crit_{\alpha}(\widehat f_{\widehat \lambda})+\delta \le \crit_{\alpha}(\widehat
f_{\lambda})+2\delta\le \crit_{\alpha}(\widehat f_{\lambda},S){+2\delta}
\]
and simple algebra that
\begin{eqnarray*}
\lefteqn{\norm{f-\Pi_{S(\widehat\lambda)}\widehat f_{\widehat \lambda}}^{2}+\alpha\norm{\widehat f_{\widehat \lambda}-\Pi_{S(\widehat\lambda)}\widehat f_{\widehat \lambda}}^{2}}\\
&\le& \norm{f-\Pi_{S}\widehat f_{\lambda}}^{2}+\alpha\norm{\widehat f_{\lambda}-\Pi_{S}\widehat f_{\lambda}}^{2}+ 2\sigma^2\penn(S) +2\delta\\
&&\ +\ 2\<\eps,\Pi_{S(\widehat\lambda)}\widehat f_{\widehat \lambda}-f\>-\sigma^2\penn(S(\widehat \lambda))\ +\ 2\<\eps, f-\Pi_{S}\widehat f_{ \lambda}\>-\sigma^2\penn(S).
\end{eqnarray*}
For $\lambda\in \Lambda$ and $S\in\SS$, let us set $u_{\lambda,S}=\pa{\Pi_{S}\widehat f_{\lambda}-f}/\norm{\Pi_{S}\widehat f_{\lambda}-f}$ if $\Pi_{S}\widehat  f_{\lambda}\ne f$ and $u_{\lambda,S}=0$ otherwise. 
For all $\lambda$ and $S$, we have $u_{\lambda,S}\in\overline S$ and
\begin{eqnarray*}
\lefteqn{\norm{f-\Pi_{S(\widehat\lambda)}\widehat f_{\widehat \lambda}}^{2}+\alpha\norm{\widehat f_{\widehat \lambda}-\Pi_{S(\widehat\lambda)}\widehat f_{\widehat \lambda}}^{2}}\\
&\le& \norm{f-\Pi_{S}\widehat f_{\lambda}}^{2}+ \alpha\norm{\widehat f_{\lambda}-\Pi_{S}\widehat f_{\lambda}}^{2}+ 2\sigma^2\penn(S){+2\delta}\\
&& +\ 2\ab{\<\eps,u_{\widehat\lambda,S(\widehat \lambda)}\>}\norm{\Pi_{S(\widehat\lambda)}\widehat f_{\widehat \lambda}-f}-\sigma^2\penn(S(\widehat \lambda))\\
&& +\ 2\ab{\<\eps,u_{\lambda,S}\>}\norm{\Pi_{S}\widehat f_{\lambda}-f}-\sigma^2\penn(S)\\
&\le& \norm{f-\Pi_{S}\widehat f_{\lambda}}^{2}+\alpha\norm{\widehat f_{\lambda}-\Pi_{S}\widehat f_{\lambda}}^{2}+2\sigma^2\penn(S){+2\delta}\\
&& +\ K^{-1}\norm{f-\Pi_{S(\widehat\lambda)}\widehat f_{\widehat \lambda}}^{2}+K\norm{\Pi_{\bar S(\widehat \lambda)}\eps}^{2}-\sigma^2\penn(S(\widehat \lambda))\\
&& +\ K^{-1}\norm{f-\Pi_{S}\widehat f_{\lambda}}^{2}+K\norm{\Pi_{\bar S}\eps}^{2}-\sigma^2\penn(S)
\end{eqnarray*}
Hence, by using~\eref{pen} and~\eref{penn} we get
\begin{eqnarray}
\lefteqn{(1-K^{-1})\norm{f-\Pi_{S(\widehat\lambda)}\widehat f_{\widehat \lambda}}^{2}+\alpha\norm{\widehat f_{\widehat \lambda}-\Pi_{S(\widehat\lambda)}\widehat f_{\widehat \lambda}}^{2}}\nonumber\\
&\le& (1+K^{-1})\norm{f-\Pi_{S}\widehat f_{\lambda}}^{2}+\alpha\norm{\widehat f_{\lambda}-\Pi_{S}\widehat f_{\lambda}}^{2}+2\sigma^2\penn(S)+\tilde  \Sigma {+2\delta}\nonumber\\
&\le& 2(1+K^{-1})\norm{f-\widehat  f_{\lambda}}^{2}{+2\delta}\nonumber\\
&& +\pa{\alpha+2(1+K^{-1})}\norm{\widehat f_{\lambda}-\Pi_{S}\widehat f_{\lambda}}^{2}+2\sigma^2\penn(S)+\tilde  \Sigma\label{eq1}
\end{eqnarray}
where 
\[
\tilde  \Sigma=2K\sum_{S\in\SS}\pa{\norm{\Pi_{\overline S}\eps}^{2}-{\ED(S)\over n-\dim(S)}\norm{Y-\Pi_{\overline S}Y}^{2}}_{+}.
\]
For each $S\in\SS$, 
\[
{\norm{Y-\Pi_{S}Y}^{2}\over n-\dim(S)}\ge {\norm{Y-\Pi_{\overline S}Y}^{2}\over n-\dim(S)}
\]
and since the variable $\norm{Y-\Pi_{\overline S}Y}^{2}$ is independent of $\norm{\Pi_{\overline S}\eps}^{2}$ and is stochastically larger than $\norm{\eps-\Pi_{\overline S}\eps}^{2}$, we deduce from the definition of $\ED(S)$ and~\eref{sigma}, that on the one hand $\E(\tilde  \Sigma)\le 2K\sigma^{2}\Sigma$. 

On the other hand, since $S$ is arbitrary among $\SS_{\lambda}$ and since 
\[
\pa{{1\over \alpha}+{1\over 1-K^{-1}}}^{-1}\norm{f-\widehat f_{\widehat \lambda}}^{2}\le (1-K^{-1})\norm{f-\Pi_{S(\widehat\lambda)}\widehat f_{\widehat \lambda}}^{2}+\alpha\norm{\widehat f_{\widehat \lambda}-\Pi_{S(\widehat\lambda)}\widehat f_{\widehat \lambda}}^{2}
\]
we deduce from~\eref{eq1} that for all $\lambda\in\Lambda$,
\begin{equation}\label{eq10}
\norm{f-\widehat f_{\widehat \lambda}}^{2}\le C^{-1}\cro{\norm{f-\widehat  f_{\lambda}}^{2}+A(\widehat f_{\lambda},\SS_{\lambda})+\tilde  \Sigma +\delta}
\end{equation}
with 
\begin{equation}
C^{-1}=C^{-1}(K,\alpha)={\pa{1+\alpha-K^{-1}}\pa{\alpha+2(1+K^{-1})}\over \alpha(1-K^{-1})},
\label{Cste.eq}
\end{equation}
and~\eref{eq3} follows by taking the expectation on both sides of~\eref{eq10}. Note that provided that 
\[
\inf_{\lambda\in\Lambda}\cro{\norm{f-\widehat  f_{\lambda}}^{2}+A(\widehat f_{\lambda},\SS_{\lambda})}
\]
is measurable, we have actually proved the stronger inequality
\begin{equation}\label{eq3bis}
C\E\cro{\norm{f-\widehat f_{\widehat \lambda}}^{2}}\le \E\cro{\inf_{\lambda\in\Lambda}\ac{\norm{f-\widehat f_{\lambda}}^{2}+A(\widehat f_{\lambda},\SS_{\lambda})}}+\sigma^{2}\Sigma +\delta.
\end{equation}

Let us now turn to the second part of the Theorem, fixing some $\lambda\in\Lambda$. Since equality holds in~\eref{pen}, under Assumption~\ref{H4} by~\eref{eq:oom} 
\[
\pen(S)=K\ED(S)\le C(\kappa,K)(\dim(S)\vee \Delta(S)),\ \ \forall S\in\SS.
\]
If $\SS_{\lambda}$ is  non-random, for some $C'=C'(\kappa,K)>0$ and all $S\in\SS_{\lambda}$,
\begin{eqnarray*}
\lefteqn{C'\E\cro{A(\widehat f_{\lambda},\SS_{\lambda})}}\\
&\le&\E\cro{\norm{\widehat f_{\lambda}-\Pi_{S}\widehat f_{\lambda}}^{2}}+(\dim(S)\vee \Delta(S))\E\cro{\widehat\sigma_{S}^{2}},\\
&=& \E\cro{\norm{\widehat f_{\lambda}-\Pi_{S}\widehat f_{\lambda}}^{2}}+{\dim(S)\vee \Delta(S)\over n-\dim(S)}\cro{\norm{f-\Pi_{S}f}^{2}+(n-\dim(S))\sigma^{2}}.
\end{eqnarray*}
Since $\norm{f-\Pi_{S}f}^{2}\le \norm{f-\Pi_{S}\widehat f_{\lambda}}^{2}$, we have
\[
\norm{f-\Pi_{S}f}^{2}\le \E\cro{\norm{f-\Pi_{S}\widehat f_{\lambda}}^{2}}\le 2\E\cro{\norm{f-\widehat f_{\lambda}}^{2}}+2\E\cro{\norm{\widehat f_{\lambda}-\Pi_{S}\widehat f_{\lambda}}^{2}},
\]
and  under Assumption~\ref{H4}, $(\dim(S)\vee \Delta(S))/(n-\dim(S))\le \kappa(1-\kappa)^{-1}$, and hence for all $S\in\SS_{\lambda}$
\begin{eqnarray*}
C'\E\cro{A(\widehat f_{\lambda},\SS_{\lambda})}&\le& \pa{1+{2\kappa\over 1-\kappa}}\E\cro{\norm{f-\widehat f_{\lambda}}^{2}}+{2\kappa\over 1-\kappa}\E\cro{\norm{\widehat f_{\lambda}-\Pi_{S}\widehat f_{\lambda}}^{2}}\\
&&\ \ +\ \ \pa{\dim(S)\vee \Delta(S)}\sigma^{2}.
\end{eqnarray*}
which leads to~\eref{BorneA1}.

Let us turn to the proof of~\eref{BorneA2}. We set $\widehat \sigma_{\lambda}^2=\widehat\sigma^2_{\widehat
S_{\lambda}}$. Since with probability one $\widehat f_{\lambda}\in\widehat S_{\lambda}\in\SS_{\lambda}$,
\[
\E\cro{A(\widehat f_{\lambda},\SS_{\lambda})}\le \E\cro{\pen(\widehat S_{\lambda})\widehat \sigma_{\lambda}^{2}}
\]
and it suffices thus to bound the right-hand side. Since equality holds in~\eref{pen} and since $\widehat f_{\lambda}\in\widehat S_{\lambda}$
\begin{eqnarray*}
\pen(\widehat S_{\lambda})\,\widehat \sigma^2_{\lambda}&=& K{\ED(\widehat S_{\lambda})\over n-\dim(\widehat S_{\lambda})}\norm{Y-\Pi_{\widehat S_{\lambda}}Y}^{2}\\
&\le&  K{\ED(\widehat S_{\lambda})\over n-\dim(\widehat S_{\lambda})}\norm{Y-\widehat f_{\lambda}}^{2}=K{\ED(\widehat S_{\lambda})\over n-\dim(\widehat S_{\lambda})}\norm{f+\eps-\widehat f_{\lambda}}^{2}\\
&\le& 2K{\ED(\widehat S_{\lambda})\over n-\dim(\widehat S_{\lambda})}\cro{\norm{f-\widehat f_{\lambda}}^{2}+\norm{\eps}^{2}}\\
&\le& 2K{\ED(\widehat S_{\lambda})\over n-\dim(\widehat S_{\lambda})}\cro{\norm{f-\widehat f_{\lambda}}^{2}+(\norm{\eps}^{2}-2n\sigma^{2})_{+}+2n\sigma^{2}}.
\end{eqnarray*}
Under Assumption~\ref{H4}, $1\le \Delta(\widehat S_{\lambda})\vee \dim(\widehat S_{\lambda})\le
\kappa n$ and we deduce from~\eref{eq:oom} that  for some constant $C$
depending only on $K$ and $\kappa$
\begin{eqnarray*}
C\pen(\widehat S_{\lambda})\, \widehat \sigma^2_{\lambda}&\le& \norm{f-\widehat f_{\lambda}}^{2}+\pa{\dim(\widehat S_{\lambda})\vee \Delta(\widehat S_{\lambda})}\sigma^{2}+(\norm{\eps}^{2}-2n\sigma^{2})_{+},
\end{eqnarray*}
and the result follows from the fact that $\E[(\norm{\eps}^{2}-2n\sigma^{2})_{+}]\le 3\sigma^{2}$ for all $n$.

\subsection{Proof of Proposition~\ref{prop-opt}}
For all $\lambda\in\Lambda$ and $f\in\ S$, $\norm{f-\widehat f_{\lambda}}\ge \norm{\Pi_{S}\widehat f_{\lambda}-\widehat f_{\lambda}}$
and hence, 
\[
\norm{f-\widehat f_{\widetilde \lambda}}^{2}\ge \inf_{\lambda\in\Lambda}\norm{f-\widehat f_{\lambda}}^{2}\ge {1\over 2}\inf_{\lambda\in\Lambda}\cro{\norm{f-\widehat f_{\lambda}}^{2}+\norm{\Pi_{S}\widehat f_{\lambda}-\widehat f_{\lambda}}^{2}}.
\]
Besides, since the minimax rate of estimation over $S$ is of order $\dim(S)\sigma^{2}$, for some universal constant $C$, 
\[
C\sup_{f\in S}\E\cro{\norm{f-\widehat f_{\widetilde \lambda}}^{2}}\ge \dim(S)\sigma^{2}.
\] 
Putting these bounds together lead to the result. 

\subsection{Proof of Proposition~\ref{prop-convex}}
Under~\eref{HM}, it is not difficult to see that $d(n,M)=n/(2\log(eM))\ge 2$ so that $\SS$ is not empty and since for all $S_{m}\in\SS_{\Cv}$
\[
(\dim(S_{m})\vee 1)\le \Delta(S_{m})=|m|+ \log\cro{\bin{M}{|m|}}\le |m|(1+\log M)\le {n\over 2},
\]
Assumptions~\ref{H0} to~\ref{H2} are satisfied with $\kappa=1/2$. Besides, the set $\Lambda_{\Cv}$ being compact, $\lambda\mapsto \crit_{\alpha}(f_{\lambda})$ admits a minimum over $\Lambda_{\Cv}$ (we shall come back the minimization of this criterion at the end of the subsection) and hence we can take $\delta=0$. By applying Theorem~\ref{main} and using~\eref{BorneA1},  the resulting estimator $\widehat f_{\Cv}=\widehat f_{\widehat \lambda}$ satisfies for some universal constant $C>0$
\begin{equation}\label{Conv1}
C\E\cro{\norm{f-\widehat f_{\Cv}}^{2}}\le \inf_{g\in\FF_{\Lambda}}\ac{\norm{f-g}^{2}+\overline A(g,\SS)},
\end{equation}
where 
\begin{equation}\label{def-Abar}
\overline A(g,\SS)=\inf_{S\in\SS}\cro{\norm{g-\Pi_{S}g}^{2}+\pa{\dim(S)\vee\Delta(S)}\sigma^{2}}.
\end{equation}
We bound $\overline A(g,\SS)$ from above by using the following approximation result below the proof of which can be found in Makovoz~\citeyearpar{MR1382053} (more precisely, we refer to the proof of his Theorem~2).
\begin{lemma}
For all $g$ in the convex hull $\FF_{\Lambda}$ of the $\phi_{j}$ and all $D\ge 1$, there exists $m\subset \ac{1,\ldots,M}$ such that $|m|=(2D)\wedge M$ and 
\[
\norm{g-\Pi_{S_{m}}g}^{2}\le 4D^{-1}\sup_{j=1,\ldots,M}\norm{\phi_{j}}^{2}.
\]
\end{lemma}
By using this lemma and the fact that $\log \bin{M}{D}\le D\log(eM/D)$ for all $D\in\{1,\ldots,M\}$, we get
\[
\overline A(g,\SS)\le \inf_{1\le D\le d(n,M)/2}\cro{{4nL^{2}\over D}+2D(1+\log(eM/(2D))}\sigma^{2}.
\]
Taking for $D$  the integer part of
\[
x(n,M,L)=\sqrt{{nL^{2}\over \log(eM/\sqrt{nL^{2}})}}
\]
which belongs to $[1,d(n,M)/2]$ under~\eref{HM}, we get
\begin{equation}\label{minConv}
\overline A(g,\SS)\le C'\sqrt{nL^{2}\log(eM/\sqrt{nL^{2}})}\sigma^{2}
\end{equation}
for some universal constant $C'>0$ which together with~\eref{Conv1} leads to the risk bound
\[
\E\cro{\norm{f-\widehat f_{\Cv}}^{2}}-C\inf_{g\in\FF_{\Lambda}}\norm{f-g}^{2}\le C\sqrt{nL^{2}\log(eM/\sqrt{nL^{2}})}\sigma^{2}.
\]

Concerning the computation of $\widehat f_{\Cv}$, note that
\begin{eqnarray*}
\inf_{\lambda\in\Lambda}\crit_{\alpha}(f_{\lambda})&=&\inf_{\lambda\in\Lambda}\inf_{S\in\SS_{\Cv}}\cro{
\norm{Y-\Pi_{S}f_{\lambda}}^{2}+\alpha\norm{f_{\lambda}-\Pi_{S} f_{\lambda}}^{2}+\pen(S)\,\widehat\sigma^2_{S}}\\
&=&\inf_{S\in\SS_{\Cv}}\ac{\cro{\inf_{\lambda\in\Lambda}\pa{
\norm{Y-\Pi_{S}f_{\lambda}}^{2}+\alpha\norm{f_{\lambda}-\Pi_{S} f_{\lambda}}^{2}}}+\pen(S)\,\widehat\sigma^2_{S}},
\end{eqnarray*}
and hence, one can solve the problem of minimizing $\crit_{\alpha}(f_{\lambda})$ over $\lambda\in\Lambda$ by proceeding into two steps. First, for each $S$ in the finite set $\SS_{\Cv}$ minimize the convex criterion 
\[
\crit_{\alpha}(S,f_{\lambda})=\norm{Y-\Pi_{S}f_{\lambda}}^{2}+\alpha\norm{f_{\lambda}-\Pi_{S} f_{\lambda}}^{2}
\]
over the convex (and compact set) $\Lambda_{\Cv}$. Denote by $\widehat f_{\Cv,S}$ the resulting minimizers. Then, minimize the quantity $\crit_{\alpha}(S,\widehat f_{\Cv,S})+\pen(S)\,\widehat\sigma^2_{S}$ for $S$ varying among $\SS_{\Cv}$. Denoting by $\widehat S$ such a minimizer, we have that $\widehat f_{\Cv}=\widehat f_{\Cv,\widehat S}$.

\subsection{Proof of Proposition~\ref{prop-simult}}
By applying Theorem~\ref{main}, we obtain that the selected estimator $\widehat f_{\widehat \lambda}$ satisfies 
\[
C\E\cro{\norm{f-\widehat f_{\widehat \lambda}}^{2}}\le \inf_{\lambda\in \ac{{\rm L},{\rm MS},\Cv} }\cro{\E\cro{\norm{f-\widehat f_{\lambda}}^{2}}+\E\cro{A(\widehat f_{\lambda},\SS_{\lambda})}}.
\]
Let us now bound $\E\cro{A(\widehat f_{\lambda},\SS_{\lambda})}$ for each $\lambda\in\Lambda$.

If $\lambda={\rm L}$, by using~\eref{BorneA1} and the fact that $\widehat f_{{\rm L}}\in S_{\{1,\ldots,M\}}$, we have
\[
C'\E\cro{A(\widehat f_{{\rm L}},\SS_{{\rm L}})}\le \E\cro{\norm{f-\widehat f_{{\rm L}}}^{2}}+M\sigma^{2}.
\]
If $\lambda={\rm MS}$,  we may use~\eref{BorneA2} since with probability one $\widehat f_{{\rm MS}}\in \SS_{{\rm MS}}$ and since $\dim(S)\vee\Delta(S)\le 1+\log(M)$ for all $S\in\SS_{{\rm MS}}$,  we get
\[
C'\E\cro{A(\widehat f_{{\rm MS}},\SS_{{\rm MS}})}\le \E\cro{\norm{f-\widehat f_{{\rm MS}}}^{2}}+\log(M)\sigma^{2}.
\]
Finally, let us turn to the case $\lambda=\Cv$ and denote by $g$ the best approximation of $f$ in ${\mathcal C}$. Since $\widehat f_{\Cv}\in {\mathcal C}$, for all $S\in\SS_{\Cv}$,
\begin{eqnarray*}
\norm{\widehat f_{\Cv}-\Pi_{S}\widehat f_{\Cv}}&\le& \norm{\widehat f_{\Cv}-\Pi_{S}g}=\norm{\widehat f_{\Cv}-f+f-g+g-\Pi_{S}g}\\
&\le& 2\norm{f-\widehat f_{\Cv}}+\norm{g-\Pi_{S}g},
\end{eqnarray*}
and hence by using~\eref{BorneA1}
\[
C'\E\cro{A(\widehat f_{\Cv},\SS_{\Cv})}\le \E\cro{\norm{f-\widehat f_{\Cv}}^{2}}+\overline A(g,\SS_{\Cv})
\]
where $\overline A(g,\SS_{\Cv})$ is given by~\eref{def-Abar}. By arguing as in Section~\eref{sect-CA}, we deduce that under~\eref{HM}
\[
C'\E\cro{A(\widehat f_{\Cv},\SS_{\Cv})}\le \E\cro{\norm{f-\widehat
   f_{\Cv}}^{2}}+\sqrt{nL^{2}\log(eM/\sqrt{nL^{2}})}\sigma^{2}.
\]
By putting these bounds together we get the result.

\subsection{Proof of Corollary~\ref{linear estimators}}
Since Assumptions~\ref{H0} to~\ref{H2} are fulfilled and $\FF$ is finite, we may apply Theorem~\ref{main} and take $\delta=0$. By using~\eref{BorneA1}, we have for some $C$ depending on $K,\alpha$ and $\kappa$,
\begin{eqnarray*}
\lefteqn{C\E\cro{\norm{f-\widehat f_{\widehat \lambda}}^{2}}}\\
&\leq& \inf_{\lambda\in\Lambda} \ac{\E\cro{\norm{f-\widehat f_{\lambda}}^2}+\E\cro{\norm{\widehat f_{\lambda}-\Pi_{S_{\lambda}}\widehat f_{\lambda}}^2}+a(1+\dim(S_{\lambda}))\sigma^2}.
\end{eqnarray*}
For all $\lambda\in\Lambda$, 
\begin{eqnarray*}
\E\cro{\|f-\widehat f_{\lambda}\|^{2}}&=&\norm{f-A_{\lambda}f}^2+\E\cro{\norm{A_{\lambda}\eps}^2}\\
&=&\norm{f-A_{\lambda}f}^2+{\rm Tr}(A_{\lambda}^*A_{\lambda})\sigma^2\\
&\ge& \max\ac{\norm{f-A_{\lambda}f}^2,{\rm Tr}(A_{\lambda}^*A_{\lambda})\sigma^2}
\end{eqnarray*}
and
\begin{eqnarray*}
\E\cro{\norm{\widehat f_{\lambda}-\Pi_{S_{\lambda}}\widehat f_{\lambda}}^2}&=&\norm{(I-\Pi_{S_{\lambda}})A_{\lambda}f}^2+\E\cro{\norm{(I-\Pi_{S_{\lambda}})A_{\lambda}\eps}^2},\\
&\le& 2\max\ac{\norm{(I-\Pi_{S_{\lambda}})A_{\lambda}f}^2,\E\cro{\norm{A_{\lambda}\eps}^2}}\\
&=& 2\max\ac{\norm{(I-\Pi_{S_{\lambda}})A_{\lambda}f}^2,{\rm Tr}(A_{\lambda}^{*}A_{\lambda})\sigma^{2}}
\end{eqnarray*}
and hence, Corollary~\ref{linear estimators} follows from the next lemma. 
\begin{lemma}\label{lemme lineaire}
For all $\lambda\in\Lambda$ we have
\begin{eqnarray*}
&(i)&\norm{(I-\Pi_{S_{\lambda}})A_{\lambda}f}\,\leq\,\norm{f-A_{\lambda}f},\\
&(ii)& \dim(S_{\lambda})\,\leq\,4\,{\rm Tr}(A_{\lambda}^*A_{\lambda}).
\end{eqnarray*}
\end{lemma}

\underline{Proof of Lemma \ref{lemme lineaire}:}
Writing
$f=f_{0}+f_{1}\in \ker(A_{\lambda})\oplus \rg(A^*_{\lambda})$
and using the fact that $ \rg(A^*_{\lambda})=\ker(A_{\lambda})^{\perp}$ and the definition of $\overline \Pi_{\lambda}$, we obtain
\begin{eqnarray*}
\norm{f-A_{\lambda}f}^2 &=& \norm{f_{0}+f_{1}-A_{\lambda}f_{1}}^2\\
&=&\norm{f_{0}-\Pi_{\ker(A_{\lambda})}A_{\lambda}f_{1}}^2+\norm{(I-\overline \Pi_{\lambda}A_{\lambda})f_{1}}^2\\
&\geq& \norm{(A_{\lambda}^+-\overline \Pi_{\lambda})A_{\lambda}f_{1}}^2\\
&\geq& \sum_{k=1}^{m_{\lambda}}s_{k}^2<A_{\lambda}f,v_{k}>^2,
\end{eqnarray*}
where $s_{1}\geq\ldots\geq s_{m_{\lambda}}$ are the singular values of $A^+_{\lambda}-\overline \Pi_{\lambda}$ counted with their multiplicity and $(v_{1},\ldots,v_{m_{\lambda}})$ is an orthonormal family of right-singular vectors associated to $(s_{1},\ldots,s_{m_{\lambda}})$. If $s_{1}<1$, then $S_{\lambda}=\R^{n}$ and we have $\norm{f-A_{\lambda}f}\ge  \norm{(I-\Pi_{S_{\lambda}})A_{\lambda}f}=0$. Otherwise, $s_{1}\ge 1$, we may consider $k_{\lambda}$ as the largest $k$ such that $s_{k}\geq 1$ and derive that
\begin{eqnarray*}
\norm{f-A_{\lambda}f}^2&\geq&\sum_{k=1}^{k_{\lambda}}s_{k}^2<A_{\lambda}f,v_{k}>^2\\
&\geq&\sum_{k=1}^{k_{\lambda}}<A_{\lambda}f,v_{k}>^2\ = \ \norm{(I-\Pi_{S_{\lambda}})A_{\lambda}f}^2,
\end{eqnarray*}
which proves the assertion~$(i)$.


For the bound~$(ii)$, 
we set $M_{\lambda}=A^+_{\lambda}-\overline \Pi_{\lambda}$ and  note that
$$(M_{\lambda}-\overline \Pi_{\lambda})(M_{\lambda}-\overline \Pi_{\lambda})^*=M_{\lambda}M^*_{\lambda}+\overline \Pi_{\lambda}\overline \Pi_{\lambda}^*-M_{\lambda}\overline \Pi_{\lambda}^*-\overline \Pi_{\lambda}M_{\lambda}^*$$
induces a  semi-positive quadratic form on $\rg(A_{\lambda}^*)$. As a consequence the quadratic form 
$(M_{\lambda}+\overline \Pi_{\lambda})(M_{\lambda}+\overline \Pi_{\lambda})^*$ is dominated by the quadratic form $2(M_{\lambda}M_{\lambda}^*+\overline \Pi_{\lambda}\overline \Pi_{\lambda}^*)$ on $\rg(A_{\lambda}^*)$.
Furthermore 
$$(M_{\lambda}+\overline \Pi_{\lambda})(M_{\lambda}+\overline \Pi_{\lambda})^*=(A_{\lambda}^+)(A_{\lambda}^+)^*=(A_{\lambda}^*A_{\lambda})^+$$
where $(A_{\lambda}^*A_{\lambda})^+$ is the inverse of  the linear operator $L_{\lambda}:\rg(A^*_{\lambda})\to\rg(A^*_{\lambda})$ induced by 
$A_{\lambda}^*A_{\lambda}$ restricted on $\rg(A^*_{\lambda})$. We then have 
that the quadratic form induced by
$(A_{\lambda}^*A_{\lambda})^+$
is dominated by the quadratic form
$$2(A_{\lambda}^{+}-\overline \Pi_{\lambda})(A_{\lambda}^{+}-\overline \Pi_{\lambda})^*+2\overline \Pi_{\lambda}\overline \Pi_{\lambda}^*$$
on $\rg(A_{\lambda}^*)$.
In particular the sequence of the eigenvalues of $(A_{\lambda}^*A_{\lambda})^+$  is dominated by the sequence $(2s_{k}^2 +2)_{k=1,m_{\lambda}}$ so
\begin{eqnarray*}
{\rm Tr}(A_{\lambda}^*A_{\lambda})\ = \ {\rm Tr}(L_{\lambda})&\geq& \sum_{k=1}^{m_{\lambda}}{1\over 2(1+s_{k}^2)}\\
&\geq&  \sum_{k=k_{\lambda}+1}^{m_{\lambda}}{1\over 2(1+s_{k}^2)}\ \geq \ {\rm dim}(S_{\lambda})/4,
\end{eqnarray*}
which conclude the proof of Lemma~\ref{lemme lineaire}.

\subsection{Proof of Corollary~\ref{cor-vs2}}
Along the section, we write $S_{*}$ for $S_{m^*}$ and $\widehat S_{\lambda}$ for $S_{\widehat m(\lambda)}$ for short.
By using~\eref{eqbase} with $\delta=0$ and  since $\Sigma\le1+\log(1+p)$, we have
$$C\E\cro{\|f-\widehat f_{\widehat \lambda}\|^2}\le \E\cro{\inf_{\lambda\in\Lambda}\|f-\Pi_{\widehat S_{\lambda}}Y\|^2+\pen(\widehat S_{\lambda})\widehat \sigma^2_{\widehat S_{\lambda}}}+(1+\log(p+1))\sigma^2,$$
for some constant $C>0$ depending on $K$ only.
Writing $B$ for the event $B=\ac{m^*\notin \widehat \M}$, we have
\[
\E\cro{\inf_{\lambda\in\Lambda}\ac{\|f-\Pi_{\widehat S_{\lambda}}Y\|^2+\pen(\widehat S_{\lambda})\widehat \sigma^2_{\widehat S_{\lambda}}}}\le A_{n}+R'_{n}
\]
where
\begin{eqnarray*}
A_{n}&=& \E\cro{\norm{f-\Pi_{S_{*}}Y}^{2}+\pen(S_{*})\widehat\sigma^2_{S_{*}}}\\
R'_{n}&=& \E\cro{\inf_{\lambda\in\Lambda}\ac{\|f-\Pi_{\widehat S_{\lambda}}Y\|^2+\pen(\widehat S_{\lambda})\widehat \sigma^2_{\widehat S_{\lambda}}}{\bf 1}_{B}}.
\end{eqnarray*}

Let us bound $A_{n}$ from above. Note that $\|f-\Pi_{S_{*}}Y\|^2=\|\Pi_{S_{*}}\eps\|^2$ and $\widehat \sigma^2_{S_{*}}=\|(I-\Pi_{S_{*}})\eps\|^2/(n-\dim(S_{*}))$ and since $\dim(S_{*})\le D_{\max}\le \kappa n/ (2\log p)$, by using~\eref{eq:oom} we get
\begin{eqnarray*}
A_{n}&\le& (\dim(S_{*})+\pen(S_{*}))\sigma^2\leq C' (1+\log(p))\dim(S_{*})\sigma^2,
\end{eqnarray*}
for some constant $C'>0$ depending on $K$ and $\kappa$ only.

Let us now turn to $R'_{n}$. For all $\lambda\in\Lambda$, $\|f-\Pi_{\widehat S_{\lambda}}Y\|^2\le \|f\|^{2}$ and 
\[
\widehat \sigma^2_{\widehat S_{\lambda}}={\|Y-\Pi_{\widehat S_{\lambda}}Y\|^{2}\over n-\dim(\widehat S_{\lambda})}\le 2{\norm{f}^{2}+\norm{\eps}^{2}\over n-\dim(\widehat S_{\lambda})}.
\]
Since for all $S\in\SS$, $\dim(S)\le D_{\max}\le  \kappa n/ (2\log p)$, by using~\eref{eq:oom} again, there exists some positive  constant $c$ depending on $K$ and $\kappa$ only such that for all $\lambda\in\Lambda$, $\pen(\widehat S_{\lambda})/(n-\dim(\widehat S_{\lambda}))\le c$ and hence, 
\[
\inf_{\lambda\in\Lambda}\ac{\|f-\Pi_{\widehat S_{\lambda}}Y\|^2+\pen(\widehat S_{\lambda})\hat \sigma^2_{\widehat S_{\lambda}}}{\bf 1}_{B}\le (1+2c)\pa{\norm{f}^{2}+\norm{\eps}^{2}}{\bf 1}_{B}.
\]
Some calculation shows that $\E\cro{\pa{\norm{f}^{2}+\norm{\eps}^{2}}^{2}}\le\pa{\norm{f}^{2}+2n\sigma^{2}}^{2}$  and hence, by Cauchy-Schwarz inequality 
\[
R'_{n}\le (1+2c)(\|f\|^2+2n\sigma^2)\sqrt{\P(B)}.
\]
The result follows by putting the bounds on $A_{n}$ and $R'_{n}$ together.


\bibliographystyle{apalike}


\section{\label{App.st}Appendix}

\subsection{\label{calcpen.st}Computation of $\ED(S)$}
The penalty  $\pen_{\Delta}(S)$,
defined at equation~\eref{EDk.eq}, is linked to the EDkhi function introduced in Baraud {\it al}~\citeyearpar{BaGiHu2009} (see Definition~3), via the following formula:
\begin{equation*}
\pen_{\Delta}(S) = {n-\dim(S)\over n-\dim(S)-1}\mathrm{EDkhi}\left( \dim(S)+1, n-\dim(S)-1, \frac{e^{-\Delta(S)}}{\dim(S)+1}\right).
\end{equation*}
Therefore, according to the result given in Section~6.1
in Baraud {\it et al}~\citeyearpar{BaGiHu2009},  $\pen_{\Delta}(S)$ is the solution in $x$ of the equation
\begin{eqnarray*}
\frac{e^{-\Delta(S)}}{D+1} &=&
\P\left(F_{D+3, N-1} \geq  x \frac{N-1}{N(D+3)}\right)
\\
&&\ \ \ \ \ -x \frac{N-1}{N(D+1)}
\P\left(F_{D+1, N+1} \geq x \frac{N+1}{N(D+1)}\right).
\end{eqnarray*}

\subsection{Simulated examples\label{SimEx.st}}

The collection $\EE$ is composed of several collections  $\EE_{1},
\ldots, \EE_{11}$ that are detailed below. The collections $\EE_{1}$ to
$\EE_{10}$ are composed of examples where $X$ is
generated as $n$ independent centered Gaussian vectors with
covariance matrix $C$. For each $e \in \{1, \ldots, 10\}$, we define
 a $p \times p$ matrix $C_{e}$ and a $p$-vector of parameters
 $\beta_{e}$. We denote by $\X_{e}$ the set of  5 matrices $X$
 simulated as $n$-i.i.d $\NN_{p}(0, C_{e})$. The collection $\EE_{e}$ is
 then defined as follows:
\begin{equation*}
\EE_{e} = \left\{ \expl(n,p,X,\beta, \rho), (n,p) \in \I, X \in
\X_{e}, \beta=\beta_{e}, \rho \in \RR \right\}
\end{equation*}
where $\RR= \{5, 10, 20\}$  and 
\begin{equation}
\I=
\left\{(100,50),(100,100),(100,1000),(200,100),(200,200)\right\}
\label{I2.eq}
\end{equation}
in Section~\ref{CVEDKhi.st}, and 
\begin{equation}
\I=
\left\{(100,50),(100,100),(200,100),(200,200)\right\} 
\label{I3.eq}
\end{equation}
in   Section~\ref{SelVar.st}.

Let us now describe the collections $\EE_{1}$ to $\EE_{10}$.
\subparagraph{ Collection $\EE_{1}$} The matrix $C$ equals the $p \times p$
identity matrix denoted $I_{p}$. 
The parameters $\beta$ satisfy
$\beta_{j}=0$ for $j \geq 16$, 
$\beta_{j}=2.5$ for $1 \leq j \leq 5$, $\beta_{j}=1.5$ for $6 \leq j
\leq 10$, $\beta_{j}=0.5$ for $11 \leq j \leq 15$.
\subparagraph{Collection $\EE_{2}$}  the  matrix $C$ is such that $C_{jk} =
r^{|j-k|}$, for $1\leq j,k \leq 15$ and $16 \leq j,k \leq p$ with
$r=0.5$. Otherwise $C_{j,k}=0$. The parameters $\beta$ are as in Collection
$\EE_{1}$.
\subparagraph{Collection $\EE_{3}$} The matrix $C$ is as in Collection
$\EE_{2}$ with $r=0.95$, the  parameters $\beta$ are as in Collection
$\EE_{1}$.
\subparagraph{Collection $\EE_{4}$} The matrix $C$ is  such that $C_{jk} =
r^{|j-k|}$, for $1\leq j,k \leq p$, with $r=0.5$, the  parameters
$\beta$ are as in Collection   $\EE_{1}$.
\subparagraph{ Collection $\EE_{5}$} the matrix $C$ is as in Collection
$\EE_{4}$ with $r=0.95$, the  parameters $\beta$ are as in Collection
$\EE_{1}$.
\subparagraph{ Collection $\EE_{6}$} The matrix $C$ equals $I_{p}$. The
parameters $\beta$ satisfy $\beta_{j}=0$ for $j \geq 16$, 
$\beta_{j}=1.5$ for $j \leq 15$.
\subparagraph{ Collection $\EE_{7}$} The matrix $C$ satisfies $C_{j,k} =
(1-\rho_{1})\1_{j=k} + \rho_{1}$ for $1 \leq, j,k \leq 3$, $C_{j,k}=C_{k,j}
= \rho_{2}$ for $j=4, k=1,2,3$, $C_{j,k}= \1_{j=k}$ for $j,k \geq
5$, with $\rho_{1}=.39$ and $\rho_{2}=.23$. The parameters 
$\beta$ satisfy $\beta_{j}=0$ for $j \geq 
4$, $\beta_{j} = 5.6$ for $j \leq 3$. 
\subparagraph{ Collection $\EE_{8}$} The matrix $C$ satisfies $C_{j,k}
=0.5^{|j-k|}$ for $j,k \leq 8$, $C_{j,k}= \1_{j=k}$ for $j,k \geq
9$. The parameters 
$\beta$ satisfy $\beta_{j}=0$ for $j \not\in \{1,2,5\}$, $\beta_{1} = 3$,
$\beta_{2}=1.5$, $\beta_{5}=2$. 
\subparagraph{ Collection $\EE_{9}$} The matrix $C$ is defined as in
Example~$\EE_{8}$. The parameters
$\beta$ satisfy $\beta_{j}=0$ for $j \geq 9$, $\beta_{j} = 0.85$ for
$j \leq 8$.
\subparagraph{ Collection $\EE_{10}$} The matrix $C$ satisfies $C_{j,k} = 0.5
\1_{j \neq k} + \1_{j=k}$
for $j,k \leq 40$, $C_{j,k} = \1_{j=k}$ for $j,k \geq 41$. The
parameters $\beta$ satisfy $\beta_{j}=2$ for $11 \leq j \leq 
20$ and $31 \leq j \leq 40$, $\beta_{j} = 0$ otherwise.
\subparagraph{ Collection $\EE_{11}$} 
In this last example, we denote
by $\X_{11}$ the set of 5 matrices $X$ simulated as follows. For $1 \leq j \leq p$, we denote by $X_{j}$ the
column $j$ of $X$. Let $E$ be generated as $n$
i.i.d. $\NN_{p}(0, 0.01 I_{p})$ and  let 
$Z_{1}, Z_{2}, Z_{3}$ be generated as $n$  i.i.d. $\NN_{3}(0,
I_{3})$. Then for $j=1, \ldots, 5$, $X_{j} = Z_{1} + E_{j}$, for
$j=6, \ldots, 10$, $X_{j} = Z_{2} + E_{j}$, for $j=11, \ldots, 15$,
$X_{j} = Z_{3} + E_{j}$, for $j \geq 16$, $X_{j}=E_{j}$.  The
parameters $\beta$ are as in Collection   $\EE_{6}$.  The collection
$\EE_{11}$ is defined as the set of
examples $\expl(n,p,X, \beta, \rho)$ for $(n,p) \in \I$, $X \in
\X_{11}$,  and $\rho \in \RR$. 

The collection $\EE$ is thus composed of 660 examples for $\I$ chosen
as in~\eref{I3.eq}, and 825  for $\I$ chosen
as in~\eref{I2.eq}. 
For some of the examples,   the Lasso estimators were highly biased
leading to high values of the ratio $O_{\expl}/n\sigma^{2}$, see
Equation~\eref{Oracl.eq}.  We only keep the examples for which the Lasso estimator improves the risk of the naive estimator $Y$ by a factor at least $1/3$. This convention leads us to remove 171 examples over 825. These pathological examples are coming from the collections
$\EE_{1}$, $\EE_{6}$ and $\EE_{7}$ for $n=100$ and $p\geq 100$, and
from collections $\EE_{2}$ and $\EE_{4}$ when $p=1000$. The examples of collection  $\EE_{7}$ were chosen by Zou to illustrate  that the Lasso estimators may be highly
biased. All the other examples, correspond to matrices $X$ that
are nearly orthogonal. 

\subsection{\label{proc.st}Procedures for calculating sets of predictors}
Let $\widehat \M=\bigcup_{\ell\in\mathcal{L}} 
\widehat \M_{\ell}$ where we recall that for $\ell\in\LL$, $\widehat \M_{\ell}=\{\widehat m(\ell,h)|\ h\in H_{\ell}\}$.

{\em The Lasso procedure} is described in
Section~\ref{CVEDKhi.st}.  The collection $\widehat \M_{{\rm Lasso}}=\{\widehat m(1),\ldots,\widehat m(D_{\max})\}$ where $\widehat{m}(h)$ is the set of indices corresponding to the predictors returned by the LARS-Lasso algorithm at step $h\in\{1,\ldots,D_{\max}\}$ (see Section~\ref{CVEDKhi.st}).

{\em The ridge procedure} is based on the minimization of $\|Y-X\beta\|^{2} + h \|\beta\|^{2}$ with respect
to $\beta$, 
for some positive $h$, see for example Hoerl and Kennard~\citeyearpar{HK06}. Tibshirani~\citeyearpar{MR1379242} noted that in
the case of a large number of small effects, ridge regression gives
better results than the lasso for variable selection. For each $h\in H_{\mathrm{ridge}}$, the regression coefficients
$\widehat{\beta}(h)$ are calculated and a  collection
of predictors sets is built as follows. Let $j_{1}, \ldots
j_{p}$ be such that $|\widehat{\beta}_{j_{1}}(h)| > \ldots >
|\widehat{\beta}_{j_{p}}(h)| $ and set  
\[
M_{h} =\ac{\{j_{1},\ldots,j_{k}\},\ k=1,\ldots,D_{\max}}.
\]
Then, the collection $\widehat\M_{\mathrm{ridge}}$ is defined as $\widehat\M_{\mathrm{ridge}} = \{
M_{h}, h \in H_{\mathrm{ridge}} \}$.

{\em The elastic net procedure} proposed by Zou and
Hastie~\citeyearpar{MR2137327}  mixes the $\ell_1$ and $\ell_2$
penalties of the Lasso and the ridge procedures. Let
$H_{\mathrm{ridge}}$ be a grid of values for the tuning parameter $h$
of the $\ell_2$ penalty. We choose
$\widehat\M_{\mathrm{en}}=\{M_{(\mathrm{en},h)} : h\in H_{\mathrm{ridge}}\}$ 
where $M_{(\mathrm{en},h)}$ denotes the collection of the active
sets of cardinality less than $D_{\max}$, selected by the elastic net
procedure when the $\ell_2$-smoothing parameter equals $h$. For each $h\in H_{\mathrm{ridge}}$ the collection  $M_{(\mathrm{en},h)}$ can be conveniently computed by first calculating the ridge regression coefficients and then applying the LARS-lasso algorithm, see Zou and Hastie~\citeyearpar{MR2137327}.

{\em The partial least squares regression} (PLSR1)
aims  to reduce the dimensionality of the regression problem by calculating
a small number of components  that are usefull for predicting
$Y$. Several applications of this procedure for analysing
high-dimensional genomic data have been reviewed by Boulesteix and
Strimmer~\citeyearpar{BS06}. In particular, it can be used  for
calculating subsets of covariates as we did for the ridge procedure.
The PLSR1
procedure
constructs, for a given $h$, 
uncorrelated latent components $t_{1}, \ldots, t_{h}$ that are
highly correlated  with the response $Y$, see
Helland~\citeyearpar{H06}.  Let
$H_{\mathrm{pls}}$ be a grid a values for the tuning parameter
$h$. For each $h \in H_{\mathrm{pls}}$, we write $\widehat \beta(h)$ for the PLS regression coefficients calculated with the first $h$ components. We then set $\widehat\M_{\mathrm{PLS}}=\{M_{h}:h\in H_{\mathrm{pls}}\}$, where  $M_{h}$
is build from $\widehat \beta(h)$ as for the ridge procedure.

{\em The adaptive lasso procedure} proposed by Zou~\citeyearpar{MR2279469} starts with
a preliminary estimator $\widetilde{\beta}$. Then one applies the
lasso procedure replacing the parameters $|\beta_{j}|, j=1, \ldots, p$
in the $\ell_1$ penalty by the weighted parameters
$|\beta_{j}|/|\tilde{\beta}_{j}|^{\gamma}, j=1, \ldots, p$ for some positive
$\gamma$. The idea is to increase the penalty for
coefficients that are close to zero, reducing thus the bias in the
estimation of $f$ and improving the variable selection accuracy. Zou showed that, if
$\widetilde{\beta}$ is a $\sqrt{n}$-consistent estimator of $\beta$,
then   the adaptive lasso procedure is consistent in situations where the
lasso is not. A lot of work has been done around this subject, see
Huang et al.~\citeyearpar{HMZ08} for example. 

We apply the procedure with $\gamma=1$, and considering two
different preliminary estimators:

\noindent
- using the ridge estimator,
$\widetilde{\beta}(h)$ as preliminary estimator. For each $h
\in H_{\mathrm{ridge}}$, the
adaptive lasso procedure is applied for calculating the active sets,
$M_{\mathrm{ALridge}, h}$, 
of cardinality less than $D_{\max}$. 
The
collection $\widehat{\M}_{\mathrm{ALridge}}$ is thus defined as 
$\widehat{\M}_{\mathrm{ALridge}}=\left\{M_{\mathrm{ALridge}, h}, h
\in H_{\mathrm{ridge}}\right\}$.

\noindent
- using the PLSR1 estimator,
$\widetilde{\beta}(h)$, as preliminary estimator. The procedure is
the same as 
described just above. The
collection $M_{\mathrm{ALpls}}$ is defined as
$M_{\mathrm{ALpls}}=\left\{M_{{\mathrm{ALpls}, h}}, h \in
H_{\mathrm{pls}}\right\}$. 

{\em The random forest algorithm} was proposed by Breiman~\citeyearpar{RandomForest} for classification and regression problems. The procedure
averages several regression trees calculated on
bootstrap samples. The algorithm returns measures of variable
importance that may be used for variable selection, see for example
D\'iaz-Uriarte and Alvares de Andr\'es~\citeyearpar{DA06}, Genuer et
al.~\citeyearpar{GePoTu}, Strobl et al.~\citeyearpar{sbzh07, sbkaz08}.

Let us denote by $h$ the number of variables randomly chosen at each split when constructing the trees and 
\[
H_{rF}= \{p/j\ |\  j \in \{3, 2, 1.5, 1\}\}.
\]
For each $h \in H_{rF}$, we consider the set of indices 
\begin{equation*}
M_h = \{\{j_1, \ldots, j_k\}, k=1,\ldots, D_{\max}\},
\end{equation*}
where $\{j_1, \ldots, j_k\}$ are the ranks of the variable importance measures. 
Two importance measures are proposed. The first one is based on the decrease in the mean
square error of prediction after permutation of each of the
variables. It leads to the collection
$\widehat{\M}_{\mathrm{rFmse}}=\{ M_h, h \in H_{rF}\}$. The second one is based on the decrease in node impurities, and leads similarly to the  collection
${\widehat{\M}_{\mathrm{purity}}}$.


{\em The exhaustive procedure} considers the collection of all subsets of
$\left\{1, \ldots p\right\}$ with dimension smaller than $D_{\max}$. We
denote this collection $\M_{\mathrm{exhaustive}}$.

\subparagraph{Choice of tuning parameters}

We have to choose  $D_{\max}$, the largest number of
predictors considered in the collection $\widehat{\M}$. For all methods, except
the exhaustive method, $D_{\max}$ may be large, say $D_{\max} \leq \min(n-2,
p)$. Nevertheless, for saving computing time, we chose $D_{\max}$ large
enough such that the dimension of the estimated subset is always
smaller than $D_{\max}$. For the exhaustive method, $D_{\max}$ must be chosen in
order to make   the calculation feasible: $D_{\max}=4$ for $p=50$, $D_{\max}=3$
for $p=100$ and $D_{\max}=2$ for $p=200$. 

For the ridge method we choose $H_{\mathrm{ridge}} = \{10^{-3}, 10^{-2},
10^{-1}, 1, 5\}$, and for the PLSR1 method, $H_{\mathrm{pls}}=1,
\ldots, 5$.

\affiliation{Universit\'e de Nice Sophia Antipolis, Ecole Polytechnique and INRA}

\address{Universit\'e de Nice Sophia-Antipolis,\\ 
Laboratoire J-A Dieudonn\'e,
UMR CNRS 6621\\
Parc Valrose\\
06108, Nice cedex 02\\
France\\
\printead{e1}\\
}

\address{Ecole Polytechnique,\\
CMAP,  UMR CNRS 7641\\
route de Saclay\\
91128 Palaiseau Cedex\\
France\\
\printead{e2}\\
}

\address{INRA MIAJ\\
78352, Jouy en Josas cedex\\
France\\
\printead{e3}\\
}

\end{document}